\documentclass[reqno]{amsart}
\usepackage{CJK}
\usepackage{hyperref}
\usepackage{cite}
\usepackage{amsmath}
\usepackage{mathrsfs}
\usepackage{xcolor}
\usepackage{bbm}
\usepackage{pifont}
\usepackage{bbding}
\usepackage{amsmath, esint}
\usepackage{mathrsfs}
\usepackage{amsfonts}
\usepackage{amssymb}
\usepackage{amsfonts,amssymb,amsmath,indentfirst,cases,amsthm}
\usepackage{epsfig}
\usepackage[numbers,sort&compress]{natbib}

\makeatletter
\@namedef{subjclassname@2020}{%
	\textup{2020} Mathematics Subject Classification}
\makeatother
\numberwithin{equation}{section}

\newtheorem{theorem}{Theorem}[section]
\newtheorem{lemma}[theorem]{Lemma}
\newtheorem{definition}[theorem]{Definition}

\newtheorem{corollary}[theorem]{Corollary}

\allowdisplaybreaks

\begin{document}
	
\title[\hfil Mixed Local and Nonlocal equations in the Heisenberg Group]{Regularity for Mixed Local and Nonlocal Degenerate Elliptic Equations in the Heisenberg Group}

\author[J. Zhang and P. Niu  \hfil \hfilneg]{Junli Zhang$^*$ and Pengcheng Niu}

\thanks{$^*$Corresponding author.}

\address{Junli Zhang \hfill\break School of Mathematics and Data Science, Shaanxi University of Science and Technology, Xi'an, Shaanxi, 710021, China}
\email{jlzhang2020@163.com}

\address{Pengcheng Niu  \hfill\break School of Mathematics and Statistics, Northwestern Polytechnical University, Xi'an, Shaanxi, 710129, China}
\email{pengchengniu@nwpu.edu.cn}

\subjclass[2020]{Primary 35B45, 35B65, 35D30; Secondary 35J92, 35R11}
\keywords{regularity, mixed local and nonlocal degenerate elliptic equation, local boundedness, H\"{o}lder continuity, Harnack inequality, weak Harnack inequality.}

\maketitle

\begin{abstract}
In this paper, we investigate the regularity for mixed local and nonlocal degenerate elliptic equations in the Heisenberg group. Inspired by the De Giorgi-Nash-Moser theory, the local boundedness of weak subsolutions and the H\"{o}lder continuity of weak solutions to mixed local and nonlocal degenerate elliptic equations are established by deriving the Caccioppoli type inequality for weak subsolutions and the logarithmic estimates for weak supersolutions. Furthermore, the Harnack inequality for weak solutions and the weak Harnack inequality for weak supersolutions are proved by using the estimates involving a Tail term and expansion of positivity.
\end{abstract}

\section{Introduction}
\label{sec-1}

For mixed local and nonlocal equations in Euclidean spaces
\begin{equation}\label{eq11}
\Delta u{\rm{ + }}P.V.\int_{{\mathbb{R}^n}} {\frac{{u(x) - u(y)}}{{{{\left| {x - y} \right|}^{n + 2s}}}}dy}  = 0,\;0 < s < 1,
\end{equation}
the Harnack inequality was proven for nonnegative solutions by Foondun in \cite{F09}, and the boundary Harnack inequality was investigated by Chen, Kim, Song and Vondra\v{c}ek in \cite{CKSV12}. For the parabolic equations corresponding to \eqref{eq11}, Barlow, Bass, Chen and Kassmann in \cite{BBC09} and Chen and Kumagai in \cite{CK10} obtained the Harnack inequality for global nonnegative solutions. Recently, Garain and Kinnunen in \cite{GK23} established a weak Harnack inequality with a remainder term for sign-changing solutions of the parabolic equations corresponding to \eqref{eq11}. Regarding the radial symmetry, maximum principle, interior and boundary Lipschitz regularity of weak solutions to \eqref{eq11}, Biagi, Dipierro, Valdinoci and Vecchi in \cite{BDVV21} and \cite{BDVV211} conducted systematic studies. For the quantitative and qualitative properties of equations with more general inhomogeneous terms $f\left( {x ,u} \right)$, as well as the research on boundary regularity for more general mixed local and nonlocal operators, one can refer to Biswas, Modasiya and Sen \cite{BMS23} and Su, Valdinoci, Wei and Zhang \cite{SVWZ22}.

For more general nonlinear mixed local and nonlocal equations in Euclidean spaces
\begin{equation}\label{eq12}
div\left( {{{\left| {\nabla u(x)} \right|}^{p - 2}}\nabla u(x)} \right)+P.V.\int_{{\mathbb{R}^n}} {\frac{{{{\left| {u(x) - u(y)} \right|}^{p - 2}}(u(x) - u(y))}}{{{{\left| {x - y} \right|}^{n + sp}}}}dy}  = 0,
\end{equation}
where $0 < s < 1,\;1 < p < \infty$, Garain and Kinnunen in \cite{GK22} proved the local boundedness, H\"{o}lder continuity, Harnack inequality and semicontinuity of weak solutions by employing the De Giorgi-Nash-Moser theory. For studies on the existence, uniqueness, and higher-order H\"{o}lder regularity, one can refer to Garain and Lindgren \cite{GL23}. Fang, Shang and Zhang in \cite{FSZ22} established the local regularity for parabolic equations corresponding to \eqref{eq12}. De Filippis and Mingione in \cite{DM} considered the functional corresponding to the mixed local and nonlocal elliptic equations with non-uniform growth
\begin{equation}\label{eq13}
a(x)div\left( {{{\left| {\nabla u(x)} \right|}^{p - 2}}\nabla u(x)} \right)+P.V.\int_{{\mathbb{R}^n}} {\frac{{{{\left| {u(x) - u(y)} \right|}^{q - 2}}(u(x) - u(y))}}{{{{\left| {x - y} \right|}^{n + sq}}}}dy}  = f\left( x \right)\;
\end{equation}
and established the local ${C^{1,\alpha }}$-regularity under the assumptions $a(x) = 1$ and $p > sq$. Ding, Fang and Zhang in \cite{DFZ24} derived the local boundedness, local H\"{o}lder continuity and Harnack inequality under the assumption $1 < p \le sq$. Byun, Lee and Song in \cite{BLS25} concluded the H\"{o}lder regularity and Harnack inequality for minimizers of the functional corresponding to \eqref{eq13} when $a(x) \ge 0$ and $f\left( x \right) = 0$.

The research on the regularity for degenerate elliptic equations formed by vector fields originated from the classic work of H\"{o}rmander \cite{H67}. There have been numerous studies on the regularity for quasilinear degenerate elliptic equations with $p$-growth conditions composed of vector fields in the Heisenberg group
\begin{equation}\label{eq10}
di{v_H}\left( {{{\left| {{\nabla _H}u(\xi )} \right|}^{p - 2}}{\nabla _H}u(\xi )} \right) = 0.
\end{equation}
Capogna in \cite{C97} first obtained the ${C^{1,\alpha }}$-regularity of weak solutions as $p=2$; Mukherjee and Zhong in \cite{MZ21} showed the ${C^{1,\alpha }}$-regularity for quasilinear degenerate elliptic equations with $p$-growth conditions as $1<p<\infty$ by using De Giorgi's method. Additionally, Zhang and Niu in \cite{ZN22} proved fractional estimates for nonlinear non-differentiable degenerate elliptic equations on the Heisenberg group. For more results, one can refer to the references listed in \cite{ZN22}. For generalized Orlicz functionals including $p,q$-growth on the Heisenberg group, Zhang and Niu in \cite{ZN20} verified the ${C^{\alpha }}$-regularity of minimizers. Zhang and Li in \cite{ZL23} established the weak differentiability of weak solutions as $1<p<4$ for quasilinear degenerate elliptic equations with $p,q$-growth conditions on the Heisenberg group; Zhang and Niu in \cite{ZN231} also considered a sufficient condition for the absence of the Lavrentiev phenomenon for quasilinear degenerate elliptic equations with $p,q$-growth conditions on the Heisenberg group.

For fractional equations in the Heisenberg group, we first note that the definition of fractional operators and the research methods are significantly different from those in Euclidean spaces. When $p=2$, Roncal and Thangavelu in \cite{RT16} defined the explicit fractional integrals in the Heisenberg group, proved several Hardy inequalities for the conformally invariant fractional sub-Laplacian and the optimality of the related constants. Some important generalized results were given by Frank, Lieb and Seiringer in \cite{FLS08}. Ferrari and Franchi in \cite{FF15} obtained Harnack and H\"{o}lder estimates for fractional equations in Carnot groups. For more results, one can also refer to \cite{FMPPS18} and \cite{FGMT15}. Manfredini, Palatucci, Piccinini and Polidoro in \cite{MPPP23} obtained the interior boundedness and H\"{o}lder continuity of fractional equations in the Heisenberg group by using the De Giorgi-Nash-Moser theory. Palatucci and Piccinini in \cite{PP22} established nonlocal Harnack inequalities in the Heisenberg group.

We cite that some extremely important tools have been introduced into the nonlocal theory in the Heisenberg group and successfully applied to fractional equations on the Heisenberg group, such as the well-known Caffarelli-Silvestre harmonic extension (see \cite{CS07}), Fourier representation, the compactness of pseudodifferential commutators in Palatucci and Pisante \cite{PP14}, and commutator estimates in Schikorra \cite{S16}, and other techniques. A special quantity, namely the nonlocal tail of a function, is introduced in these references, which plays a crucial role in overcoming the nonlocality of nonlinear operators.

Palatucci and Piccinini in 2022 (see \cite{PP22}) and 2023 (see \cite{PP23}) proposed an open problem how to determine the regularity of nonlocal double-phase equations in the Heisenberg group. Fang, Zhang and Zhang provided a solution in \cite{FZZ24}. For the regularity to mixed local and nonlocal degenerate elliptic equations in the Heisenberg group
\begin{equation}\label{eq14}
a\left( \xi  \right)di{v_H}\left( {{{\left| {{\nabla _H}u(\xi )} \right|}^{p - 2}}{\nabla _H}u(\xi )} \right)+b\left( \xi  \right)P.V.\int_{{{\rm{H}}^n}} {\frac{{|u(\xi ) - u(\eta ){|^{q - 2}}(u(\xi ) - u(\eta ))}}{{\left\| {{\eta ^{ - 1}} \circ \xi } \right\|_{{{\rm{H}}^n}}^{Q + sq}}}d\eta }  = f\left( \xi  \right),
\end{equation}
Zhang, Niu and Wu in \cite{ZNW25} proved that both the derivative in the vertical direction and the derivative in the horizontal direction of weak solutions $u \in H{W^{1,p}}\left( \Omega  \right)$ to \eqref{eq14} with $p=q=2$ and $f \in {L^2}\left( \Omega  \right)$ belonging to the first-order local Sobolev space $HW_{loc}^{1,p}\left( \Omega  \right)$. Furthermore, under the assumption $f \in H{W^{m,p}}\left( \Omega  \right)$, we can also demonstrate $u \in HW_{loc}^{m + 2,p}\left( \Omega  \right)$.

In this paper, we are concerned with the regularity to the mixed local and nonlocal degenerate elliptic equation \eqref{eq14} with $p = q,\;a\left( \xi  \right) = b\left( \xi  \right) = 1$, i.e.,
\begin{equation}\label{eq15}
 - di{v_H}\left( {{{\left| {{\nabla _H}u\left( \xi  \right)} \right|}^{p - 2}}{\nabla _H}u\left( \xi  \right)} \right) + {\rm{P}}{\rm{.V}}{\rm{.}}\int_{{{\rm{H}}^n}} {\frac{{{{\left| {u\left( \xi  \right) - u\left( \eta  \right)} \right|}^{p - 2}}\left( {u\left( \xi  \right) - u\left( \eta  \right)} \right)}}{{\left\| {{\eta ^{ - 1}} \circ \xi } \right\|_{{{\rm{H}}^n}}^{Q + sp}}}d\eta }  = 0,\;\xi  \in \Omega ,
\end{equation}
where $\Omega $ is a bounded domain in ${{\rm{H}}^n}$, $1 < p < \infty ,$ $0 < s < 1,$ $Q = 2n + 2$ denotes the homogeneous dimension of ${{\rm{H}}^n}$, and P.V. stands for the Cauchy principal value.

The main results are as follows:

\begin{theorem}[Local boundedness of weak subsolutions]\label{Th11}
Let $u \in HW_{loc}^{1,p}\left( \Omega  \right)$ $(1 < p < \infty )$ be a weak subsolution to \eqref{eq15}, and ${B_r} \equiv {B_r}\left( {{\xi _0}} \right) \subset \Omega ,\;r \in \left( {0,1} \right],$ here ${B_r}$ is a ball defined by the C-C metric. Then for $\delta  \in \left( {0,1} \right]$, there exists a positive constant $c = c\left( {n,p,s} \right)$ such that
\begin{equation}\label{eq16}
\mathop {{\rm{ess}}\sup u}\limits_{{B_{\frac{r}{2}}}\left( {{\xi _0}} \right)}  \le \delta {\rm{Tail}}({u_ + };{\xi _0},\frac{r}{2}) + c{\delta ^{ - \frac{{\left( {p - 1} \right)\kappa }}{{p\left( {\kappa  - 1} \right)}}}}{\left( {\fint_{{B_r}} {u_ + ^pd\xi } } \right)^{\frac{1}{p}}},
\end{equation}
where
\[\kappa  = \left\{ \begin{array}{l}
\frac{Q}{{Q - p}},\;\;\;\;\;\;\;\;1 < p < Q,\\
\frac{{qQ}}{{p\left( {Q - q} \right)}},\;\;\;p \ge Q,\;\frac{{Qp}}{{Q + p}} < q < Q,
\end{array} \right.\;\;\;\]
and ${\rm{Tail}}\left(  \cdot  \right)$ is given by \eqref{eq21} below.
\end{theorem}

\begin{theorem}[H\"{o}lder continuity]\label{Th12}
Let $u \in HW_{loc}^{1,p}\left( \Omega  \right)$ $\left( {1 < p < Q} \right)$ be a weak solution to \eqref{eq15}. Then $u$ is locally H\"{o}lder continuous in $\Omega $, that is, for ${B_{2r}} \equiv {B_{2r}}\left( {{\xi _0}} \right) \subset \Omega ,\;r \in \left( {0,1} \right],\;\rho  \in \left( {0,r} \right]$, there exist constants $\alpha  \in \left( {0,\frac{{sp}}{{p - 1}}} \right)$ and $c = c\left( {n,p,s} \right)$ such that
\begin{equation}\label{eq17}
\mathop {{\rm{osc}}u}\limits_{{B_\rho }\left( {{\xi _0}} \right)}  = \mathop {{\rm{ess}}\sup }\limits_{{B_\rho }\left( {{\xi _0}} \right)} u - \mathop {{\rm{ess inf}}}\limits_{{B_\rho }\left( {{\xi _0}} \right)} u \le c{\left( {\frac{\rho }{r}} \right)^\alpha }\left[ {{\rm{Tail}}(u;{\xi _0},r) + {{\left( {\fint_{{B_{2r}}} {{{\left| u \right|}^p}d\xi } } \right)}^{\frac{1}{p}}}} \right],
\end{equation}
where ${\rm{Tail}}\left(  \cdot  \right)$ is given by \eqref{eq21}.
\end{theorem}

\begin{theorem}[Harnack inequality]\label{Th13}
Let $u \in HW_{loc}^{1,p}\left( \Omega  \right)$ $\left( {1 < p < Q} \right)$ be a weak solution to \eqref{eq15} satisfying $u \ge 0$ in ${B_R}\left( {{\xi _0}} \right) \subset \Omega $. Then there exists a positive constant $c = c\left( {n,p,s} \right)$ such that
\begin{equation}\label{eq18}
\mathop {{\rm{ess sup}}}\limits_{{B_{\frac{r}{2}}}\left( {{\xi _0}} \right)} u \le c\mathop {{\rm{ess inf}}}\limits_{{B_r}\left( {{\xi _0}} \right)} u + c{\left( {\frac{r}{R}} \right)^{\frac{p}{{p - 1}}}}{\rm{Tail}}({u_ - };{\xi _0},R),
\end{equation}
where $r \in \left( {0,1} \right]$ and ${B_r}\left( {{\xi _0}} \right) \subset {B_{\frac{R}{2}}}\left( {{\xi _0}} \right)$. Here ${\rm{Tail}}\left(  \cdot  \right)$ is given by \eqref{eq21}.
\end{theorem}

\begin{theorem}[Weak Harnack inequality]\label{Th14}
Let $u \in HW_{loc}^{1,p}\left( \Omega  \right)$ $\left( {1 < p < Q} \right)$ be a weak supersolution to \eqref{eq15} satisfying $u \ge 0$ in ${B_R}\left( {{\xi _0}} \right) \subset \Omega $. Then there exists a positive constant $c = c\left( {n,p,s} \right)$ such that
\begin{equation}\label{eq19}
{\left( {\fint_{{B_{\frac{r}{2}}}\left( {{\xi _0}} \right)} {{u^l}d\xi } } \right)^{\frac{1}{l}}} \le c\mathop {{\rm{ess inf}}}\limits_{{B_r}\left( {{\xi _0}} \right)} u + c{\left( {\frac{r}{R}} \right)^{\frac{p}{{p - 1}}}}{\rm{Tail}}({u_ - };{\xi _0},R),
\end{equation}
where $r \in \left( {0,1} \right]$, ${B_r}\left( {{\xi _0}} \right) \subset {B_{\frac{R}{2}}}\left( {{\xi _0}} \right)$ and $0 < l < \frac{{p*\left( {p - 1} \right)}}{p},\;p*=\frac{pQ}{Q-p}$. Here ${\rm{Tail}}\left(  \cdot  \right)$ is given by \eqref{eq21}.
\end{theorem}

The method adopted in this paper is based on the De Giorgi-Nash-Moser theory. Theorem \ref{Th11} is proved by using the Sobolev embedding theorem (Lemma \ref{Le22}), the Caccioppoli-type inequality (Lemma \ref{Le27}), and the iteration lemma (Lemma \ref{Le26}). It is worth mentioning that although the Sobolev embedding theorem only holds for $1 < p < Q$, when $p\ge Q$, we can still establish the local boundedness of weak subsolutions (Theorem \ref{Th11}) by choosing $\frac{{Qp}}{{Q + p}} < q < Q$, then applying the Sobolev embedding theorem to exponent $q$ and combining it with H\"{o}lder's inequality. That is, Theorem \ref{Th11} holds for all $1 < p < \infty $. Then, by virtue of the logarithmic lemma (Lemma \ref{Le41}), the Poincar\'{e}-type inequality (Lemma \ref{Le21}), Theorem \ref{Th11}, Lemma \ref{Le22} and Lemma \ref{Le27}, an oscillation estimate (Lemma \ref{Le43}) is derived from which Theorem \ref{Th12} directly follows. Finally, the Harnack inequality for weak solutions and the weak Harnack inequality for weak supersolutions are proven by using the estimates involving a Tail term, expansion of positivity, Lemma \ref{Le21}, Theorem \ref{Th11} and Lemma \ref{Le22}. Since the Poincar\'{e}-type inequality and the Sobolev embedding theorem only hold for $1 < p < Q$, Theorems \ref{Th12}-\ref{Th14} are only valid for $1 < p < Q$.

The paper is organized as follows: In Section 2, we introduce relevant knowledge of the Heisenberg group, function spaces, and several necessary lemmas. In particular, a Caccioppoli-type inequality for weak subsolutions to \eqref{eq15} is proven. In Section 3, we present the proof of Theorem \ref{Th11}. In Section 4, we prove Theorem \ref{Th12}. The proofs of Theorem \ref{Th13} and Theorem \ref{Th14} are given in Sections 5 and 6,  respectively.

\section{Preliminaries}
\label{Section 2}

In this section, we introduce briefly some relevant knowledge of the Heisenberg group ${{\rm{H}}^n}$, some function spaces and useful lemmas and then prove a Caccioppoli-type inequality for weak subsolutions to \eqref{eq15}.

The Euclidean space ${{\mathbb{R}}^{2n + 1}}\;(n \ge 1)$ with the group multiplication
\[\;\xi \circ \eta = \left( {{x_1} + {y_1},{x_2} + {y_2}, \cdots ,{x_{2n}} + {y_{2n}},\tau + \tau' + \frac{1}{2}\sum\limits_{i = 1}^n {\left( {{x_i}{y_{n + i}} - {x_{n + i}}{y_i}} \right)} } \right),\]
where $\xi = \left( {{x_1},{x_2}, \cdots ,{x_{2n}},\tau} \right),$ $\eta = \left( {{y_1},{y_2}, \cdots ,{y_{2n}},\tau'} \right) \in {{\mathbb{R}}^{2n+1}},$ leads to the Heisenberg group ${{\rm{H}}^n}$. The left invariant vector field on ${{\rm{H}}^n}$ is of the form
\begin{equation*}
  {X_i} = {\partial _{{x_i}}} - \frac{{{x_{n + i}}}}{2}{\partial _\tau},\;{X_{n + i}} = {\partial _{{x_{n + i}}}} + \frac{{{x_i}}}{2}{\partial _\tau},\quad 1 \le i \le n
\end{equation*}
and a non-trivial commutator is
\begin{equation*}
  T = {\partial _\tau} = \left[ {{X_i},{X_{n + i}}} \right] = {X_i}{X_{n + i}} - {X_{n + i}}{X_i},~1 \le i \le n.
\end{equation*}
We call that ${X_1},{X_2}, \cdots ,{X_{2n}}$ are the horizontal vector fields on ${{\mathbb{H}}^n}$ and $T$ the vertical vector field. Denote the horizontal gradient of a smooth function $u$ on ${{\rm{H}}^n}$ by
\begin{equation*}
{\nabla _H}u = \left( {{X_1}u,{X_2}u, \cdots ,{X_{2n}}u} \right).
\end{equation*}

The Haar measure in ${{\rm{H}}^n}$ is equivalent to the Lebesgue measure in ${{\mathbb{R}}^{2n+1}}$. We denote the Lebesgue measure of a measurable set $E \subset {{\rm{H}}^n}$ by $\left| E \right|$. For $\xi = \left( {{x_1},{x_2}, \cdots ,{x_{2n}},\tau} \right),$ we define its module as
\[{\| \xi \|_{{{\rm{H}}^n}}} = {\left( {{{\left( {\sum\limits_{i = 1}^{2n} {{x_i}^2} } \right)}^2} + {\tau^2}} \right)^{\frac{1}{4}}}.\]
The Carnot-Carath\'{e}odary metric between two points $\xi$ and $\eta$ in ${{\rm{H}}^n}$ is the shortest length of the horizontal curve joining them, denoted by $d(\xi,\eta)$. The C-C metric is equivalent to the Kor\`{a}nyi metric, i.e., $d\left( {\xi,\eta} \right) \sim {\| {{\xi^{ - 1}}\circ \eta} \|_{{{\rm{H}}^n}}}$. The ball
\begin{equation*}
  {B_r }\left( \xi_0 \right) = \left\{ {\xi \in {{\rm{H}}^n}:d\left( {\xi,\xi_0} \right) < r } \right\}
\end{equation*}
is defined by the C-C metric $d$. When not important or clear from the context, we will omit the center as follows: $B_r:=B_r( \xi_0)$. The homogeneous dimension of ${{\rm{H}}^n}$ is $Q = 2n + 2$.

Let $1 \le p < \infty ,\;\Omega  \subset {\rm{H}^n}.$ The Sobolev space $H{W^{k,p}}\left( \Omega  \right)$ is defined by
\[H{W^{k,p}}\left( \Omega  \right) = \left\{ {u \in {L^p}\left( \Omega  \right):{\nabla _H}u \in {L^p}\left( \Omega  \right),\nabla _H^2u \in {L^p}\left( \Omega  \right), \cdots ,\nabla _H^ku \in {L^p}\left( \Omega  \right)} \right\}\]
endowed with the norm
\[{\left\| u \right\|_{H{W^{k,p}}\left( \Omega  \right)}} = {\left\| u \right\|_{{L^p}\left( \Omega  \right)}} + \sum\limits_{m = 1}^k {{{\left\| {\nabla _H^mu} \right\|}_{{L^p}\left( \Omega  \right)}}} .\]
This makes $H{W^{k,p}}\left( \Omega  \right)$ a Banach space. Let $HW_0^{k,p}\left( \Omega  \right)$ denote the closure of $C_0^\infty \left( \Omega  \right)$ in $H{W^{k,p}}\left( \Omega  \right)$. The local Sobolev space $HW_{loc}^{k,p}\left( \Omega  \right)$ is defined as
\[HW_{loc}^{k,p}\left( \Omega  \right): = \left\{ {u:u \in H{W^{k,p}}\left( {\Omega '} \right),\forall \Omega ' \subset  \subset \Omega } \right\}.\]

Let $1 \le p < \infty ,\;s \in \left( {0,1} \right)$, and $v:{{\rm{H}}^n} \to {\mathbb{R}}$ be a measurable function. The Gagliardo semi-norm of $v$ is defined as
\[{\left[ v \right]_{H{W^{s,p}}\left( {{{\mathbb{H}}^n}} \right)}} = {\left( {\int_{{{\mathbb{H}}^n}} {\int_{{{\mathbb{H}}^n}} {\frac{{{{\left| {v\left( \xi  \right) - v\left( \eta  \right)} \right|}^p}}}{{\| {{\eta ^{ - 1}} \circ \xi } \|_{{{\mathbb{H}}^n}}^{Q + sp}}}\,d\xi } d\eta } } \right)^{\frac{1}{p}}},\]
and the fractional Sobolev spaces $H{W^{s,p}}\left( {{{\mathbb{H}}^n}} \right)$ on the Heisenberg group are defined as
\[H{W^{s,p}}\left( {{{\mathbb{H}}^n}} \right) = \left\{ {v \in {L^p}\left( {{{\mathbb{H}}^n}} \right):{{\left[ v \right]}_{H{W^{s,p}}\left( {{{\mathbb{H}}^n}} \right)}} < \infty } \right\},\]
endowed with the natural fractional norm
\[{\| v \|_{H{W^{s,p}}\left( {{{\mathbb{H}}^n}} \right)}} = {\left( {\| v \|_{{L^p}\left( {{{\mathbb{H}}^n}} \right)}^p + \left[ v \right]_{H{W^{s,p}}\left( {{{\mathbb{H}}^n}} \right)}^p} \right)^{\frac{1}{p}}}.\]
For any open set $\Omega  \subset {{\mathbb{H}}^n}$, we can define similarly fractional Sobolev spaces $H{W^{s,p}}\left( \Omega  \right)$ and fractional norm ${\| v \|_{H{W^{s,p}}\left( \Omega  \right)}}$. The space $HW_0^{s,p}\left( \Omega  \right)$ is the closure of $C_0^\infty \left( \Omega  \right)$ in $H{W^{s,p}}\left( \Omega  \right)$. Throughout this paper, we denote a generic positive constant as $c$ or $C$. If necessary, relevant dependencies on parameters will be illustrated by parentheses, i.e., $c=c(n,p)$ means that $c$ depends on $n,p$.

For any $u \in H{W^{s,p}}\left( {{{\rm{H}}^n}} \right)$ and any ${B_r}\left( {{\xi _0}} \right) \subset {{\rm{H}}^n}$, the nonlocal tail of the function $u$ with respect to a ball ${B_r}\left( {{\xi _0}} \right)$ is defined by
\begin{equation}\label{eq21}
{\rm{Tail}}\left( {u,{\xi _0},r} \right) = {\left( {{r^p}{\int _{{{\rm{H}}^n}\backslash {B_r}\left( {{\xi _0}} \right)}}\frac{{{{\left| {u\left( \xi  \right)} \right|}^{p - 1}}}}{{\left\| {\xi _0^{ - 1} \circ \xi } \right\|_{{{\rm{H}}^n}}^{Q + sp}}}d\xi } \right)^{\frac{1}{{p - 1}}}}.
\end{equation}

\begin{lemma}[Poincar\'{e}-type inequality, (1.1) in \cite{DLS07}]\label{Le21}
Let $1 \le p < Q$ and $v \in H{W^{1,p}}({B_r})$, then it follows
\[{\left( {\int_{{B_r}} {{{\left| {v - {{\left( v \right)}_r}} \right|}^{p*}}{\mkern 1mu} d\xi } } \right)^{\frac{1}{{p*}}}} \le c{\left( {\int_{{B_r}} {{{\left| {{\nabla _H}v} \right|}^p}d\xi } } \right)^{\frac{1}{p}}},\]
where  $p* = \frac{{pQ}}{{Q - p}},\;\;c = c(n,p) > 0,$ ${\left( v \right)_r} = \fint_{{B_r}} {vd\xi } .$
\end{lemma}

\begin{lemma}[Sobolev inequality, \cite{CDG93}]\label{Le22}
Let $1 < p < Q$ and ${B_r} \subset {{\rm{H}}^n}$. Then for any $u \in HW_0^{1,p}\left( {{B_r}} \right),$ there exists a positive constant $c = c\left( {n,p} \right)$ such that
\[{\left( {\fint_{{B_r}} {{{\left| u \right|}^{p*}}d\xi } } \right)^{\frac{1}{{p*}}}} \le cr{\left( {\fint_{{B_r}} {{{\left| {{\nabla _H}u} \right|}^p}d\xi } } \right)^{\frac{1}{p}}},\]
where $p* = \frac{{pQ}}{{Q - p}}$.
\end{lemma}

\begin{lemma}[\cite{MPPP23}]\label{Le23}
Let $1 < p < \infty ,\;0 < {s_1} \le s < 1$, $\Omega  \subset {{\rm{H}}^n}$ and $u \in H{W^{s,p}}\left( \Omega  \right)$. Then, for a positive constant $c$ depending only on $n,p,{s_1}$, it holds
\[{\left\| u \right\|_{H{W^{{s_1},p}}\left( \Omega  \right)}} \le c{\left\| u \right\|_{H{W^{s,p}}\left( \Omega  \right)}}.\]
\end{lemma}

\begin{definition}\label{De24}
Let $u \in HW_{loc}^{1,p}\left( \Omega  \right)$. The function $u$ is called a weak subsolution (supersolution) to \eqref{eq15} if for any $\Omega ' \subset  \subset \Omega $ and non-negative test function $\phi  \in HW_0^{1,p}\left( {\Omega '} \right)$, the inequality
\begin{align}\label{eq22}
   & \int_{\Omega '} {{{\left| {{\nabla _H}u} \right|}^{p - 2}}{\nabla _H}u \cdot {\nabla _H}\phi d\xi } \nonumber \\
  + & \int_{{{\rm{H}}^n}} {\int_{{{\rm{H}}^n}} {\frac{{{{\left| {u\left( \xi  \right) - u\left( \eta  \right)} \right|}^{p - 2}}\left( {u\left( \xi  \right) - u\left( \eta  \right)} \right)\left( {\phi \left( \xi  \right) - \phi \left( \eta  \right)} \right)}}{{\left\| {{\eta ^{ - 1}} \circ \xi } \right\|_{{{\rm{H}}^n}}^{Q + sp}}}d\xi d\eta } }  \le \left(  \ge  \right)0
\end{align}
holds. If $u$ is both a weak subsolution and a weak supersolution to \eqref{eq15}, then $u$ is a weak solution to \eqref{eq15}.
\end{definition}

It is obvious that if $u$ is a weak subsolution (supersolution) to \eqref{eq15}, then $-u$ is a weak supersolution (subsolution) to \eqref{eq15}.

\begin{lemma}[\cite{GK22}]\label{Le25}
Let $u$ be a weak subsolution to \eqref{eq15}, then ${u_ + } = \max \{ u,0\} $ is also a weak subsolution to \eqref{eq15}.
\end{lemma}

\begin{lemma}[\cite{D93}]\label{Le26}
Let $\left( {{Y_j}} \right)_{j = 0}^\infty $ be a sequence of positive real numbers. If there exist constants ${c_0},b > 1$ and $\beta  > 0$ such that
\begin{center}
${Y_0} \le c_0^{ - \frac{1}{\beta }}{b^{ - \frac{1}{{{\beta ^2}}}}}$ and ${Y_{j + 1}} \le {c_0}{b^j}Y_j^{1 + \beta },\;j = 0,1,2, \cdots $,
\end{center}
then
\[\mathop {\lim }\limits_{j \to \infty } {Y_j} = 0.\]
\end{lemma}

Next, we will prove

\begin{lemma}[Caccioppoli-type inequality]\label{Le27}
Let $u \in HW_{loc}^{1,p}\left( \Omega  \right)\;\;(1 < p < \infty )$ be a weak subsolution to \eqref{eq15}, and $\omega  = {\left( {u - k} \right)_ + },\;k \in \mathbb{R}$. Then for any ${B_r} \equiv {B_r}\left( {{\xi _0}} \right) \subset \Omega $ and non-negative function $\psi  \in C_0^\infty \left( {{B_r}} \right)$, it holds
\begin{align}\label{eq23}
   & \int_{{B_r}} {{\psi ^p}{{\left| {{\nabla _H}\omega } \right|}^p}d\xi }  + \int_{{B_r}} {\int_{{B_r}} {\frac{{{{\left| {\omega \left( \xi  \right)\psi \left( \xi  \right) - \omega \left( \eta  \right)\psi \left( \eta  \right)} \right|}^p}}}{{\left\| {{\eta ^{ - 1}} \circ \xi } \right\|_{{{\rm{H}}^n}}^{Q + sp}}}{\mkern 1mu} d\xi d\eta } } \nonumber \\
 \le &  c(\int_{{B_r}} {{\omega ^p}{{\left| {{\nabla _H}\psi } \right|}^p}d\xi }  + \int_{{B_r}} {\int_{{B_r}} {\max {{\{ \omega \left( \xi  \right),\omega \left( \eta  \right)\} }^p}\frac{{{{\left| {\psi \left( \xi  \right) - \psi \left( \eta  \right)} \right|}^p}}}{{\left\| {{\eta ^{ - 1}} \circ \xi } \right\|_{{{\rm{H}}^n}}^{Q + sp}}}{\mkern 1mu} d\xi d\eta } } \nonumber\\
 & + \mathop {{\rm{ess}}\;\sup }\limits_{\xi  \in {\rm{supp}}\;\psi } \int_{{{\rm{H}}^n}\backslash {B_r}} {\frac{{\omega {{\left( \eta  \right)}^{p - 1}}}}{{\left\| {{\eta ^{ - 1}} \circ \xi } \right\|_{{{\rm{H}}^n}}^{Q + sp}}}d\eta }  \cdot \int_{{B_r}} {\omega {\psi ^p}d\xi } ),
\end{align}
where $c = c\left( p \right)$. If $u$ is a weak supersolution to \eqref{eq15}, then \eqref{eq23} holds for $\omega  = {\left( {u - k} \right)_ - },\;k \in \mathbb{R}$.
\end{lemma}

\begin{proof}
Let $u \in HW_{loc}^{1,p}\left( \Omega  \right)$ be a weak subsolution to \eqref{eq15} and $\psi  \in C_0^\infty \left( {{B_r}} \right)$ be a non-negative function. For $\omega  = {\left( {u - k} \right)_ + },\;k \in \mathbb{R}$, we take $\phi  = \omega {\psi ^p}$ as the test function in \eqref{eq22}, then it follows
\begin{align}\label{eq24}
   0 \ge& \int_{{B_r}} {{{\left| {{\nabla _H}u} \right|}^{p - 2}}{\nabla _H}u \cdot {\nabla _H}\left( {\omega {\psi ^p}} \right)d\xi } \nonumber \\
   &  + \int_{{{\rm{H}}^n}} {\int_{{{\rm{H}}^n}} {\frac{{{{\left| {u\left( \xi  \right) - u\left( \eta  \right)} \right|}^{p - 2}}\left( {u\left( \xi  \right) - u\left( \eta  \right)} \right)\left( {\omega \left( \xi  \right)\psi {{\left( \xi  \right)}^p} - \omega \left( \eta  \right)\psi {{\left( \eta  \right)}^p}} \right)}}{{\left\| {{\eta ^{ - 1}} \circ \xi } \right\|_{{{\rm{H}}^n}}^{Q + sp}}}d\xi d\eta } }\nonumber \\
  = &I + J.
\end{align}

For $I$, we obtain by using Young's inequality that
\begin{align}\label{eq25}
  I & = \int_{{B_r}} {{{\left| {{\nabla _H}u} \right|}^{p - 2}}{\nabla _H}u \cdot {\nabla _H}\left( {\omega {\psi ^p}} \right)d\xi }  \nonumber  \\
   &  = \int_{{B_r}} {{{\left| {{\nabla _H}u} \right|}^{p - 2}}{\nabla _H}u \cdot {\nabla _H}\omega {\psi ^p}d\xi }  + p\int_{{B_r}} {{{\left| {{\nabla _H}u} \right|}^{p - 2}}{\nabla _H}u \cdot \omega {\psi ^{p - 1}}{\nabla _H}\psi d\xi } \nonumber \\
   &\ge \int_{{B_r}} {{{\left| {{\nabla _H}\omega } \right|}^p}{\psi ^p}d\xi }  - p\int_{{B_r}} {{{\left| {{\nabla _H}\omega } \right|}^{p - 1}}{\psi ^{p - 1}} \cdot \omega \left| {{\nabla _H}\psi } \right|d\xi }  \nonumber \\
   & \ge c\left( p \right)\int_{{B_r}} {{{\left| {{\nabla _H}\omega } \right|}^p}{\psi ^p}d\xi }  - \int_{{B_r}} {{\omega ^p}{{\left| {{\nabla _H}\psi } \right|}^p}d\xi } .
\end{align}
For $J$, it follows from the proof of Theorem 1.3 in \cite{MPPP23} that
\begin{align}\label{eq26}
   J =& \int_{{{\rm{H}}^n}} {\int_{{{\rm{H}}^n}} {\frac{{{{\left| {u\left( \xi  \right) - u\left( \eta  \right)} \right|}^{p - 2}}\left( {u\left( \xi  \right) - u\left( \eta  \right)} \right)\left( {\omega \left( \xi  \right)\psi {{\left( \xi  \right)}^p} - \omega \left( \eta  \right)\psi {{\left( \eta  \right)}^p}} \right)}}{{\left\| {{\eta ^{ - 1}} \circ \xi } \right\|_{{{\rm{H}}^n}}^{Q + sp}}}d\xi d\eta } }  \nonumber\\
   \ge&  {c_1}\left( p \right)\int_{{B_r}} {\int_{{B_r}} {\frac{{{{\left| {\omega \left( \xi  \right)\psi \left( \xi  \right) - \omega \left( \eta  \right)\psi \left( \eta  \right)} \right|}^p}}}{{\left\| {{\eta ^{ - 1}} \circ \xi } \right\|_{{{\rm{H}}^n}}^{Q + sp}}}{\mkern 1mu} d\xi d\eta } } \nonumber\\
   & - {c_2}\left( p \right)\int_{{B_r}} {\int_{{B_r}} {\max {{\{ \omega \left( \xi  \right),\omega \left( \eta  \right)\} }^p}\frac{{{{\left| {\psi \left( \xi  \right) - \psi \left( \eta  \right)} \right|}^p}}}{{\left\| {{\eta ^{ - 1}} \circ \xi } \right\|_{{{\rm{H}}^n}}^{Q + sp}}}{\mkern 1mu} d\xi d\eta } } \nonumber\\
   & - 2\int_{{{\rm{H}}^n}\backslash {B_r}} {\int_{{B_r}} {\frac{{\omega {{\left( \eta  \right)}^{p - 1}}\omega \left( \xi  \right)\psi {{\left( \xi  \right)}^p}}}{{\left\| {{\eta ^{ - 1}} \circ \xi } \right\|_{{{\rm{H}}^n}}^{Q + sp}}}d\xi d\eta } } .
\end{align}

Substituting \eqref{eq25} and \eqref{eq26} into \eqref{eq24} yields that \eqref{eq23} holds. By applying Lemma \ref{Le27} to $-u$, we obtain the result for weak supersolutions.
\end{proof}

\section{Local Boundedness of Weak Subsolutions}
\label{Section 3}

In this section, we first establish an iterative scheme by using Lemma \ref{Le22} and Lemma \ref{Le27}, and then prove Theorem \ref{Th11} by virtue of Lemma \ref{Le26}. Since Lemma \ref{Le22} only holds for $1 < p < Q$, the proof will be carried out in two cases: $1 < p < Q$ and $p \ge Q$.

\textbf{Proof of Theorem \ref{Th11}.}
For $j = 0,1,2, \cdots $, denote
\begin{center}
${r_j} = \frac{r}{2}\left( {1 + {2^{ - j}}} \right),\;{\bar r_j} = \frac{{{r_j} + {r_{j + 1}}}}{2},\;{B_j} = {B_{{r_j}}}\left( {{\xi _0}} \right)$ and ${\bar B_j} = {B_{{{\bar r}_j}}}\left( {{\xi _0}} \right)$.
\end{center}
Let $\left( {{\psi _j}} \right)_{j = 0}^\infty  \subset C_0^\infty \left( {{{\bar B}_j}} \right)\;\;(j = 0,1,2, \cdots )$ be a sequence of cut-off functions satisfying $0 \le {\psi _j} \le 1,\;\left| {{\nabla _H}{\psi _j}} \right| \le \frac{{{2^{j + 3}}}}{r}$ in ${\bar B_j}$ and ${\psi _j} = 1$ in ${B_{j + 1}}$. Moreover, for $k,\bar k \ge 0$, denote
\begin{center}
${k_j} = k + \left( {1 - {2^{ - j}}} \right)\bar k,\;{\bar k_j} = \frac{{{k_j} + {k_{j + 1}}}}{2},\;{\omega _j} = {\left( {u - {k_j}} \right)_ + }$ and ${\bar \omega _j} = {\left( {u - {{\bar k}_j}} \right)_ + }.$
\end{center}
Next, we prove Theorem \ref{Th11} by dividing it into two cases: $1 < p < Q$ and $p \ge Q$.

\textbf{Case 1.} When $1 < p < Q$, on the one hand, applying Lemma \ref{Le22} to the function ${\bar \omega _j}{\psi _j}$ yields
\begin{align*}
   & {\left( {\fint_{{B_j}} {{{\left| {{{\bar \omega }_j}{\psi _j}} \right|}^{p*}}d\xi } } \right)^{\frac{p}{{p*}}}} \le {cr_j^p}\fint_{{B_j}} {{{\left| {{\nabla _H}\left( {{{\bar \omega }_j}{\psi _j}} \right)} \right|}^p}d\xi }  \\
  = &{cr_j^p} \fint_{{B_j}} {\psi _j^p{{\left| {{\nabla _H}{{\bar \omega }_j}} \right|}^p}dx}  + {cr_j^p}\fint_{{B_j}} {\bar \omega _j^p{{\left| {{\nabla _H}{\psi _j}} \right|}^p}dx} .
\end{align*}
Combining with Lemma \ref{Le27}, we obtain
\begin{align}\label{eq31}
  {\left( {\fint_{{B_j}} {{{\left| {{{\bar \omega }_j}{\psi _j}} \right|}^{p*}}d\xi } } \right)^{\frac{p}{{p*}}}} \le  & {cr_j^p}\fint_{{B_j}} {{{\bar \omega }_j}^p{{\left| {{\nabla _H}{\psi _j}} \right|}^p}d\xi }  \nonumber\\
   &  +{cr_j^p} \int_{{B_j}} {\fint_{{B_j}} {\max {{\{ {{\bar \omega }_j}\left( \xi  \right),{{\bar \omega }_j}\left( \eta  \right)\} }^p}\frac{{{{\left| {{\psi _j}\left( \xi  \right) - {\psi _j}\left( \eta  \right)} \right|}^p}}}{{\left\| {{\eta ^{ - 1}} \circ \xi } \right\|_{{{\rm{H}}^n}}^{Q + sp}}}{\mkern 1mu} d\xi d\eta } }\nonumber\\
   & +{cr_j^p} \mathop {{\rm{ess}}\;\sup }\limits_{\xi  \in {\rm{supp}}\;{\psi _j}} \int_{{{\rm{H}}^n}\backslash {B_j}} {\frac{{{{\bar \omega }_j}{{\left( \eta  \right)}^{p - 1}}}}{{\left\| {{\eta ^{ - 1}} \circ \xi } \right\|_{{{\rm{H}}^n}}^{Q + sp}}}d\eta }  \cdot \fint_{{B_j}} {{{\bar \omega }_j}\psi _j^pd\xi } \nonumber\\
    = :&{I_1} + {I_2} + {I_3}.
\end{align}
For ${I_1}$, by using the properties of ${\psi _j}$, $\frac{r}{2} \le {r_j} \le r$ and ${{\bar \omega }_j} \le {{ \omega }_j}$, we obtain
\begin{equation}\label{eq32}
{I_1} \le {c{2^{pj}}}\fint_{{B_j}} {\omega _j^pd\xi } .
\end{equation}
By the proof of Theorem 1.1 in \cite{MPPP23}, we have
\begin{equation}\label{eq33}
{I_2} \le c{2^{pj}}\fint_{{B_j}} {\omega _j^p\left( \xi  \right)d\xi }
\end{equation}
and
\begin{equation}\label{eq34}
{I_3} \le \frac{{c{2^{j\left( {Q + sp + p - 1} \right)}}}}{{{{\bar k}^{p - 1}}}}{\left[ {{\rm{Tail}}\left( {{\omega _0};{\xi _0},\frac{r}{2}} \right)} \right]^{p - 1}}\fint_{{B_j}} {\omega _j^p\left( \xi  \right)d\xi } .
\end{equation}
Substituting \eqref{eq32}-\eqref{eq34} into \eqref{eq31}, we get
\begin{equation}\label{eq35}
{\left( {\fint_{{B_j}} {{{\left| {{{\bar \omega }_j}{\psi _j}} \right|}^{p*}}d\xi } } \right)^{\frac{p}{{p*}}}} \le c{2^{j\left( {Q + sp + p - 1} \right)}}\left( {1 + \frac{{{{\left[ {\rm{Tail}\left( {{\omega _0};{\xi _0},\frac{r}{2}} \right)} \right]}^{p - 1}}}}{{{{\bar k}^{p - 1}}}}} \right)\fint_{{B_j}} {\omega _j^p\left( \xi  \right)d\xi } .
\end{equation}

On the other hand, using that ${\psi _j} = 1$ on ${B_{j + 1}}$, ${\bar \omega _j} \ge {\omega _{j + 1}}$ and ${\bar \omega _j} \ge {k_{j + 1}} - {\bar k_j}$ as ${\omega _{j + 1}} \ge 0$, we have
\begin{align}\label{eq36}
   & {\left( {\fint_{{B_j}} {{{\left| {{{\bar \omega }_j}{\psi _j}} \right|}^{p*}}d\xi } } \right)^{\frac{p}{{p*}}}} \ge c{\left( {\fint_{{B_{j + 1}}} {{{\left| {{{\bar \omega }_j}{\psi _j}} \right|}^{p*}}d\xi } } \right)^{\frac{p}{{p*}}}} \nonumber\\
   =&  c{\left( {\fint_{{B_{j + 1}}} {{{\left| {{{\bar \omega }_j}} \right|}^{p* - p}}{{\left| {{{\bar \omega }_j}} \right|}^p}d\xi } } \right)^{\frac{p}{{p*}}}} \ge c{\left( {\fint_{{B_{j + 1}}} {{{\left| {{{\bar \omega }_j}} \right|}^{p* - p}}\omega _{j + 1}^p\left( \xi  \right)d\xi } } \right)^{\frac{p}{{p*}}}}\nonumber\\
    \ge& c{\left( {{k_{j + 1}} - {{\bar k}_j}} \right)^{\frac{{\left( {p* - p} \right)p}}{{p*}}}}{\left( {\fint_{{B_{j + 1}}} {\omega _{j + 1}^p\left( \xi  \right)d\xi } } \right)^{\frac{p}{{p*}}}}\nonumber\\
 = &c{\left( {\frac{{\bar k}}{{{2^{j + 2}}}}} \right)^{\frac{{\left( {p* - p} \right)p}}{{p*}}}}{\left( {\fint_{{B_{j + 1}}} {\omega _{j + 1}^p\left( \xi  \right)d\xi } } \right)^{\frac{p}{{p*}}}}.
\end{align}
Therefore, combining \eqref{eq35} and \eqref{eq36}, we gain
\begin{align}\label{eq37}
   & {\left( {\frac{{\bar k}}{{{2^{j + 2}}}}} \right)^{\frac{{\left( {p* - p} \right)p}}{{p*}}}}{\left( {\fint_{{B_{j + 1}}} {\omega _{j + 1}^p\left( \xi  \right)d\xi } } \right)^{\frac{p}{{p*}}}} \nonumber\\
   \le &  c{2^{j\left( {Q + sp + p - 1} \right)}}\left( {1 + \frac{{{{\left[ {\rm{Tail}\left( {{\omega _0};{\xi _0},\frac{r}{2}} \right)} \right]}^{p - 1}}}}{{{{\bar k}^{p - 1}}}}} \right)\fint_{{B_j}} {\omega _j^p\left( \xi  \right)d\xi } .
\end{align}
Denoting ${A_j} = {\left( {\fint_{{B_j}} {\omega _j^p\left( \xi  \right)d\xi } } \right)^{\frac{1}{p}}}$, the above formula becomes
\begin{equation}\label{eq38}
{\left( {\frac{{{{\bar k}^{1 - \frac{p}{{p*}}}}}}{{{2^{\left( {j + 2} \right)\frac{{p* - p}}{{p*}}}}}}} \right)^p}A_{j + 1}^{\frac{{{p^2}}}{{p*}}} \le c{2^{j\left( {Q + sp + p - 1} \right)}}\left( {1 + \frac{{{{\left[ {{\rm{Tail}}\left( {{\omega _0};{\xi _0},\frac{r}{2}} \right)} \right]}^{p - 1}}}}{{{{\bar k}^{p - 1}}}}} \right)A_j^p.
\end{equation}

Taking
\begin{equation}\label{eq39}
\bar k \ge \delta {\rm{Tail}}\left( {{\omega _0};{\xi _0},\frac{r}{2}} \right),\;\delta  \in \left( {0,1} \right],
\end{equation}
we have
\begin{equation}\label{eq310}
{\left( {\frac{{{A_{j + 1}}}}{{\bar k}}} \right)^{\frac{p}{{p*}}}} \le {\delta ^{\frac{{1 - p}}{p}}}{\bar c^{\frac{p}{{p*}}}}{2^{j\left( {\frac{p}{Q} + \frac{{Q + sp + p - 1}}{p}} \right)}}\frac{{{A_j}}}{{\bar k}},
\end{equation}
where $\bar c = {c^{\frac{{p*}}{{{p^2}}}}}{\left[ {{\delta ^{p - 1}} + 1} \right]^{\frac{{p*}}{{{p^2}}}}}{2^{\frac{{2\left( {p* - p} \right)}}{{p*}}\frac{{p*}}{p}}}.$ Taking $C: = {2^{\frac{p}{{Q - p}} + \frac{{Q\left( {Q + sp + p - 1} \right)}}{{p\left( {Q - p} \right)}}}} > 1,\;\beta  = \frac{{p*}}{p} - 1,$ the formula \eqref{eq310} becomes
\begin{equation}\label{eq311}
\frac{{{A_{j + 1}}}}{{\bar k}} \le {\delta ^{\frac{{p*\left( {1 - p} \right)}}{{{p^2}}}}}\bar c{C^j}{\left( {\frac{{{A_j}}}{{\bar k}}} \right)^{1 + \beta }}.
\end{equation}
Moreover, since
\[\frac{{p*\left( {p - 1} \right)}}{{{p^2}\beta }} = \frac{{\left( {p - 1} \right)Q}}{{{p^2}}},\]
and in consideration of \eqref{eq39}, we choose
\begin{equation}\label{eq312}
\bar k = \delta {\rm{Tail}}\left( {{\omega _0};{\xi _0},\frac{r}{2}} \right) + {\delta ^{ - \frac{{\left( {p - 1} \right)Q}}{{{p^2}}}}}H{A_0},\;\;H = {\bar c^{\frac{1}{\beta }}}{C^{\frac{1}{{{\beta ^2}}}}}
\end{equation}
to see
\begin{equation}\label{eq313}
\frac{{{A_0}}}{{\bar k}} \le {\delta ^{\frac{{p*\left( {p - 1} \right)}}{{{p^2}\beta }}}}{\bar c^{ - \frac{1}{\beta }}}{C^{ - \frac{1}{{{\beta ^2}}}}}.
\end{equation}
Hence, by applying Lemma \ref{Le26}, we have
\[\mathop {\lim }\limits_{j \to \infty } {A_j} = {\left( {\fint_{{B_{\frac{r}{2}}}} {\left( {u\left( \xi  \right) - k - \bar k} \right)_ + ^pd\xi } } \right)^{\frac{1}{p}}} = 0,\]
i.e.
\begin{align*}
   & \mathop {\sup }\limits_{{B_{\frac{r}{2}}}} u\left( \xi  \right) \le k + \bar k \\
  = &  k + \delta {\rm{Tail}}\left( {{{\left( {u - k} \right)}_ + };{\xi _0},\frac{r}{2}} \right) + {\delta ^{ - \frac{{\left( {p - 1} \right)Q}}{{{p^2}}}}}H{\left( {\fint_{{B_r}} {\left( {u\left( \xi  \right) - k} \right)_ + ^pd\xi } } \right)^{\frac{1}{p}}}.
\end{align*}
If we take $k=0$, then \eqref{eq16} holds.

\textbf{Case 2.} When $p \ge Q$, we choose $\frac{{Qp}}{{Q + p}} < q < Q$. On the one hand, applying Lemma \ref{Le22} and H\"{o}lder's inequality to the function ${\bar \omega _j}{\psi _j}$ yields
\begin{align*}
  {\left( {\fint_{{B_j}} {{{\left| {{{\bar \omega }_j}{\psi _j}} \right|}^{q*}}d\xi } } \right)^{\frac{1}{{q*}}}} & \le c{r_j}{\left( {\fint_{{B_j}} {{{\left| {{\nabla _H}\left( {{{\bar \omega }_j}{\psi _j}} \right)} \right|}^q}d\xi } } \right)^{\frac{1}{q}}} \\
   &  \le c{r_j}{\left( {\fint_{{B_j}} {{{\left| {{\nabla _H}\left( {{{\bar \omega }_j}{\psi _j}} \right)} \right|}^p}d\xi } } \right)^{\frac{1}{p}}},
\end{align*}
so
\begin{align*}
   {\left( {\fint_{{B_j}} {{{\left| {{{\bar \omega }_j}{\psi _j}} \right|}^{q*}}d\xi } } \right)^{\frac{p}{{q*}}}}& \le {cr_j^p}\fint_{{B_j}} {{{\left| {{\nabla _H}\left( {{{\bar \omega }_j}{\psi _j}} \right)} \right|}^p}d\xi }  \\
   &  = {cr_j^p}\fint_{{B_j}} {\psi _j^p{{\left| {{\nabla _H}{{\bar \omega }_j}} \right|}^p}dx}  + {cr_j^p}\fint_{{B_j}} {\bar \omega _j^p{{\left| {{\nabla _H}{\psi _j}} \right|}^p}dx} .
\end{align*}
Combining Lemma \ref{Le27} with \eqref{eq31}-\eqref{eq34}, it follows
\begin{equation}\label{eq314}
{\left( {\fint_{{B_j}} {{{\left| {{{\bar \omega }_j}{\psi _j}} \right|}^{q*}}d\xi } } \right)^{\frac{p}{{q*}}}} \le c{2^{j\left( {Q + sp + p - 1} \right)}}\left( {1 + \frac{{{{\left[ {\rm{Tail}\left( {{\omega _0};{\xi _0},\frac{r}{2}} \right)} \right]}^{p - 1}}}}{{{{\bar k}^{p - 1}}}}} \right)\fint_{{B_j}} {\omega _j^p\left( \xi  \right)d\xi } .
\end{equation}

On the other hand, similar to the estimates in \eqref{eq36} and noting $\frac{{Qp}}{{Q + p}} < q$ (i.e. $q* > p$), we immediately deduce
\begin{equation}\label{eq315}
{\left( {\fint_{{B_j}} {{{\left| {{{\bar \omega }_j}{\psi _j}} \right|}^{q*}}d\xi } } \right)^{\frac{p}{{q*}}}} \ge c{\left( {\frac{{\bar k}}{{{2^{j + 2}}}}} \right)^{\frac{{\left( {q* - p} \right)p}}{{q*}}}}{\left( {\fint_{{B_{j + 1}}} {\omega _{j + 1}^p\left( \xi  \right)d\xi } } \right)^{\frac{p}{{q*}}}}.
\end{equation}
Combining \eqref{eq314} and \eqref{eq315}, we get
\begin{align}\label{eq316}
   & {\left( {\frac{{\bar k}}{{{2^{j + 2}}}}} \right)^{\frac{{\left( {q* - p} \right)p}}{{q*}}}}{\left( {\fint_{{B_{j + 1}}} {\omega _{j + 1}^p\left( \xi  \right)d\xi } } \right)^{\frac{p}{{q*}}}} \nonumber\\
   \le&  c{2^{j\left( {Q + sp + p - 1} \right)}}\left( {1 + \frac{{{{\left[ {\rm{Tail}\left( {{\omega _0};{\xi _0},\frac{r}{2}} \right)} \right]}^{p - 1}}}}{{{{\bar k}^{p - 1}}}}} \right)\fint_{{B_j}} {\omega _j^p\left( \xi  \right)d\xi } .
\end{align}

Similar to the derivation of \eqref{eq310}, if we take
\begin{equation}\label{eq317}
\bar k = \delta {\rm{Tail}}\left( {{\omega _0};{\xi _0},\frac{r}{2}} \right) + {\delta ^{ - \frac{{q*\left( {p - 1} \right)}}{{p\left( {q* - p} \right)}}}}H{A_0},\;\;H = {\bar c^{\frac{1}{\beta }}}{C^{\frac{1}{{{\beta ^2}}}}},\;\delta  \in \left( {0,1} \right],
\end{equation}
then
\begin{equation}\label{eq318}
{\left( {\frac{{{A_{j + 1}}}}{{\bar k}}} \right)^{\frac{p}{{q*}}}} \le {\delta ^{\frac{{1 - p}}{p}}}{\bar c^{\frac{p}{{q*}}}}{2^{j\left( {\frac{{q* - p}}{{q*}} + \frac{{Q + sp + p - 1}}{p}} \right)}}\frac{{{A_j}}}{{\bar k}},
\end{equation}
where $\bar c = {c^{\frac{{q*}}{{{p^2}}}}}{\left[ {{\delta ^{p - 1}} + 1} \right]^{\frac{{q*}}{{{p^2}}}}}{2^{\frac{{2\left( {q* - p} \right)}}{{q*}}\frac{{q*}}{p}}}.$ If we take $C: = {2^{\left( {\frac{{q* - p}}{{q*}} + \frac{{Q + sp + p - 1}}{p}} \right)\frac{{q*}}{p}}} > 1,\;\beta  = \frac{{q*}}{p} - 1,$ then the formula \eqref{eq318} becomes
\begin{equation}\label{eq319}
\frac{{{A_{j + 1}}}}{{\bar k}} \le {\delta ^{\frac{{q*\left( {1 - p} \right)}}{{{p^2}}}}}\bar c{C^j}{\left( {\frac{{{A_j}}}{{\bar k}}} \right)^{1 + \beta }}.
\end{equation}
Moreover, since
\[\frac{{q*\left( {p - 1} \right)}}{{{p^2}\beta }} = \frac{{q*\left( {p - 1} \right)}}{{p\left( {q* - p} \right)}},\]
and in consideration of \eqref{eq317}, we have
\begin{equation}\label{eq320}
\frac{{{A_0}}}{{\bar k}} \le {\delta ^{\frac{{q*\left( {p - 1} \right)}}{{{p^2}\beta }}}}{\bar c^{ - \frac{1}{\beta }}}{C^{ - \frac{1}{{{\beta ^2}}}}}.
\end{equation}
Therefore, by applying Lemma \ref{Le26}, we obtain
\[\mathop {\lim }\limits_{j \to \infty } {A_j} = {\left( {\fint_{{B_{\frac{r}{2}}}} {\left( {u\left( \xi  \right) - k - \bar k} \right)_ + ^pd\xi } } \right)^{\frac{1}{p}}} = 0,\]
i.e.
\begin{align*}
   & \mathop {\sup }\limits_{{B_{\frac{r}{2}}}} u\left( \xi  \right) \le k + \bar k \\
  = &  k + \delta {\rm{Tail}}\left( {{{\left( {u - k} \right)}_ + };{\xi _0},\frac{r}{2}} \right) + {\delta ^{ - \frac{{q*\left( {p - 1} \right)}}{{p\left( {q* - p} \right)}}}}H{\left( {\fint_{{B_r}} {\left( {u\left( \xi  \right) - k} \right)_ + ^pd\xi } } \right)^{\frac{1}{p}}}.
\end{align*}
If we take $k=0$, then \eqref{eq16} holds.

Now combining Case 1 and Case 2, we see that Theorem \ref{Th11} holds.

\section{H\"{o}lder Continuity}
\label{Section 4}

In this section, we prove Theorem \ref{Th12}. Since Theorem \ref{Th12} can be directly derived from the oscillation estimate (Lemma \ref{Le43}), the key of this section is to prove Lemma \ref{Le43}. Before presenting the proof of Lemma \ref{Le43}, we first give a logarithmic-type lemma (Lemma \ref{Le41}). Then, by using this lemma and Lemma \ref{Le21}, we obtain Corollary \ref{Co42}. Finally, we prove Lemma \ref{Le43} by applying the induction, iteration techniques, Theorem \ref{Th11}, Corollary \ref{Co42}, Lemma \ref{Le27} and Lemma \ref{Le22}.

\begin{lemma}[Logarithmic-type lemma]\label{Le41}
Let $u \in HW_{loc}^{1,p}\left( \Omega  \right)\;\;(1 < p < \infty )$ be a weak supersolution to \eqref{eq15} and satisfy $u \ge 0$ in ${B_R}\left( {{\xi _0}} \right) \subset \Omega $. Then, for any ${B_r} \equiv {B_r}\left( {{\xi _0}} \right) \subset {B_{\frac{R}{2}}}\left( {{\xi _0}} \right),\;r \in \left( {0,1} \right]$ and $d > 0$, there exists a positive constant $c = c\left( {n,p,s} \right)$ such that
\begin{align}\label{eq41}
   & \int_{{B_r}} {{{\left| {{\nabla _H}\log \left( {u + d} \right)} \right|}^p}d\xi }  + \int_{{B_r}} {\int_{{B_r}} {{{\left| {\log \left( {\frac{{u\left( \xi  \right) + d}}{{u\left( \eta  \right) + d}}} \right)} \right|}^p}\frac{1}{{\left\| {{\eta ^{ - 1}} \circ \xi } \right\|_{{{\rm{H}}^n}}^{Q + sp}}}d\xi d\eta } } \nonumber \\
  \le &  c{r^Q}\left( {{r^{ - p}} + {d^{1 - p}}{R^{ - p}}{\rm{Tail}}{{({u_ - };{\xi _0},R)}^{p - 1}}} \right).
\end{align}
Here ${\rm{Tail}}\left(  \cdot  \right)$ is given by \eqref{eq21}.
\end{lemma}

\begin{proof}
Let $\psi  \in C_0^\infty \left( {{B_{\frac{{3r}}{2}}}\left( {{\xi _0}} \right)} \right)$ satisfy
\begin{center}
$0 \le \psi  \le 1,\;\left| {{\nabla _H}\psi } \right| \le \frac{8}{r}$ in ${B_{\frac{{3r}}{2}}}\left( {{\xi _0}} \right)$ and $\psi  = 1$ in ${B_r}\left( {{\xi _0}} \right)$.
\end{center}
Choosing $\phi  = {\left( {u + d} \right)^{1 - p}}{\psi ^p}$ as the test function in \eqref{eq22}, we have
\begin{align}\label{eq42}
   0 \le& \int_{{B_{2r}}} {{{\left| {{\nabla _H}u} \right|}^{p - 2}}{\nabla _H}u \cdot {\nabla _H}\left( {{{\left( {u + d} \right)}^{1 - p}}{\psi ^p}} \right)d\xi } \nonumber \\
  &  + \int_{{{\rm{H}}^n}} {\int_{{{\rm{H}}^n}} {\frac{{{{\left| {u\left( \xi  \right) - u\left( \eta  \right)} \right|}^{p - 2}}\left( {u\left( \xi  \right) - u\left( \eta  \right)} \right)\left( {{{\left( {u\left( \xi  \right) + d} \right)}^{1 - p}}\psi {{\left( \xi  \right)}^p} - {{\left( {u\left( \eta  \right) + d} \right)}^{1 - p}}\psi {{\left( \eta  \right)}^p}} \right)}}{{\left\| {{\eta ^{ - 1}} \circ \xi } \right\|_{{{\rm{H}}^n}}^{Q + sp}}}d\xi d\eta } } \nonumber \\
 =&\int_{{B_{2r}}} {{{\left| {{\nabla _H}\left( {u + d} \right)} \right|}^{p - 2}}{\nabla _H}\left( {u + d} \right) \cdot {\nabla _H}\left( {{{\left( {u + d} \right)}^{1 - p}}{\psi ^p}} \right)d\xi }  \nonumber \\
  & + \int_{{B_{2r}}} {\int_{{B_{2r}}} {\frac{{{{\left| {u\left( \xi  \right) - u\left( \eta  \right)} \right|}^{p - 2}}\left( {u\left( \xi  \right) - u\left( \eta  \right)} \right)\left( {{{\left( {u\left( \xi  \right) + d} \right)}^{1 - p}}\psi {{\left( \xi  \right)}^p} - {{\left( {u\left( \eta  \right) + d} \right)}^{1 - p}}\psi {{\left( \eta  \right)}^p}} \right)}}{{\left\| {{\eta ^{ - 1}} \circ \xi } \right\|_{{{\rm{H}}^n}}^{Q + sp}}}d\xi d\eta } } \nonumber \\
  &  + 2\int_{{{\rm{H}}^n}\backslash {B_{2r}}} {\int_{{B_{2r}}} {\frac{{{{\left| {u\left( \xi  \right) - u\left( \eta  \right)} \right|}^{p - 2}}\left( {u\left( \xi  \right) - u\left( \eta  \right)} \right){{\left( {u\left( \xi  \right) + d} \right)}^{1 - p}}\psi {{\left( \xi  \right)}^p}}}{{\left\| {{\eta ^{ - 1}} \circ \xi } \right\|_{{{\rm{H}}^n}}^{Q + sp}}}d\xi d\eta } } \nonumber \\
   =& {I_1} + {I_2} + I{}_3.
\end{align}

For term ${I_1}$, by Young's inequality, it gets
\begin{align}\label{eq43}
   {I_1} =& \int_{{B_{2r}}} {{{\left| {{\nabla _H}\left( {u + d} \right)} \right|}^{p - 2}}{\nabla _H}\left( {u + d} \right) \cdot {\nabla _H}\left( {{{\left( {u + d} \right)}^{1 - p}}{\psi ^p}} \right)d\xi } \nonumber \\
  = &  \left( {1 - p} \right)\int_{{B_{2r}}} {{{\left| {{\nabla _H}\left( {u + d} \right)} \right|}^p}{{\left( {u + d} \right)}^{ - p}}{\psi ^p}d\xi }\nonumber \\
  & + p\int_{{B_{2r}}} {{{\left| {{\nabla _H}\left( {u + d} \right)} \right|}^{p - 2}}{\nabla _H}\left( {u + d} \right) \cdot {{\left( {u + d} \right)}^{1 - p}}{\psi ^{p - 1}}{\nabla _H}\psi d\xi } \nonumber \\
   \le& \left( {1 - p} \right)\int_{{B_{2r}}} {{{\left( {u + d} \right)}^{ - p}}{{\left| {{\nabla _H}\left( {u + d} \right)} \right|}^p}{\psi ^p}d\xi } \nonumber \\
 & + \frac{{p - 1}}{2}\int_{{B_{2r}}} {{{\left( {u + d} \right)}^{ - p}}{{\left| {{\nabla _H}\left( {u + d} \right)} \right|}^p}{\psi ^p}d\xi }  + c\left( p \right)\int_{{B_{2r}}} {{{\left| {{\nabla _H}\psi } \right|}^p}d\xi }  \nonumber \\
  \le&  - \frac{{p - 1}}{2}\int_{{B_{2r}}} {{{\left| {{\nabla _H}\log \left( {u + d} \right)} \right|}^p}{\psi ^p}d\xi }  + c\left( p \right){r^{Q - p}}.
\end{align}

For ${I_2}$ and ${I_3}$, it follows from the estimates of ${I_1}$ and ${I_2}$ in the proof of Lemma 1.4 in \cite{MPPP23} that
\begin{align}\label{eq44}
   {I_2}& \le  - \frac{1}{{c\left( p \right)}}\int_{{B_{2r}}} {\int_{{B_{2r}}} {{{\left| {\log \left( {\frac{{u\left( \xi  \right) + d}}{{u\left( \eta  \right) + d}}} \right)} \right|}^p}\frac{{\psi {{\left( \eta  \right)}^p}}}{{\left\| {{\eta ^{ - 1}} \circ \xi } \right\|_{{{\rm{H}}^n}}^{Q + sp}}}d\xi d\eta } }  + c\left( p \right){r^{Q - sp}} \nonumber\\
   &  \le  - \frac{1}{{c\left( p \right)}}\int_{{B_{2r}}} {\int_{{B_{2r}}} {{{\left| {\log \left( {\frac{{u\left( \xi  \right) + d}}{{u\left( \eta  \right) + d}}} \right)} \right|}^p}\frac{{\psi {{\left( \eta  \right)}^p}}}{{\left\| {{\eta ^{ - 1}} \circ \xi } \right\|_{{{\rm{H}}^n}}^{Q + sp}}}d\xi d\eta } }  + c\left( p \right){r^{Q - p}}
\end{align}
and
\begin{align}\label{eq45}
   {I_3}& \le c{r^{Q - sp}} + c{d^{1 - p}}{r^Q}{R^{ - sp}} + c{d^{1 - p}}{r^Q}{R^{ - p}}{\rm{Tail}}{({u_ - };{\xi _0},R)^{p - 1}}\nonumber \\
   &  \le c{r^{Q - p}} + c{d^{1 - p}}{r^Q}{R^{ - p}}{\rm{Tail}}{({u_ - };{\xi _0},R)^{p - 1}}.
\end{align}

Substituting \eqref{eq43}-\eqref{eq45} into \eqref{eq42} and using $\psi  = 1$ on ${B_r}\left( {{\xi _0}} \right)$, we have
\begin{align*}
   & \int_{{B_r}} {{{\left| {{\nabla _H}\log \left( {u + d} \right)} \right|}^p}d\xi }  + \int_{{B_r}} {\int_{{B_r}} {{{\left| {\log \left( {\frac{{u\left( \xi  \right) + d}}{{u\left( \eta  \right) + d}}} \right)} \right|}^p}\frac{1}{{\left\| {{\eta ^{ - 1}} \circ \xi } \right\|_{{{\rm{H}}^n}}^{Q + sp}}}d\xi d\eta } }  \\
   =&  \int_{{B_r}} {{{\left| {{\nabla _H}\log \left( {u + d} \right)} \right|}^p}{\psi ^p}d\xi }  + \int_{{B_r}} {\int_{{B_r}} {{{\left| {\log \left( {\frac{{u\left( \xi  \right) + d}}{{u\left( \eta  \right) + d}}} \right)} \right|}^p}\frac{{\psi {{\left( \eta  \right)}^p}}}{{\left\| {{\eta ^{ - 1}} \circ \xi } \right\|_{{{\rm{H}}^n}}^{Q + sp}}}d\xi d\eta } } \\
    \le &\int_{{B_{2r}}} {{{\left| {{\nabla _H}\log \left( {u + d} \right)} \right|}^p}{\psi ^p}d\xi }  + \int_{{B_{2r}}} {\int_{{B_{2r}}} {{{\left| {\log \left( {\frac{{u\left( \xi  \right) + d}}{{u\left( \eta  \right) + d}}} \right)} \right|}^p}\frac{{\psi {{\left( \eta  \right)}^p}}}{{\left\| {{\eta ^{ - 1}} \circ \xi } \right\|_{{{\rm{H}}^n}}^{Q + sp}}}d\xi d\eta } } \\
     \le& c{r^{Q - p}} + c{d^{1 - p}}{r^Q}{R^{ - p}}{\rm{Tail}}{({u_ - };{\xi _0},R)^{p - 1}}\\
      = &c{r^Q}\left( {{r^{ - p}} + {d^{1 - p}}{R^{ - p}}{\rm{Tail}}{{({u_ - };{\xi _0},R)}^{p - 1}}} \right),
\end{align*}
then \eqref{eq41} holds.
\end{proof}

Similar to the proof of Corollary 4.1 in \cite{MPPP23}, we use the Poincar\'{e} type inequality (Lemma \ref{Le21}) and Lemma \ref{Le41} to obtain the following corollary.

\begin{corollary}\label{Co42}
Let $u \in HW_{loc}^{1,p}\left( \Omega  \right)\;\;(1 < p < Q)$ be a weak supersolution to \eqref{eq15} and satisfy $u \ge 0$ in ${B_R}\left( {{\xi _0}} \right) \subset \Omega $. For $a,d > 0,\;b > 1$, denote
\[v = \min \left\{ {{{\left( {\log \left( {\frac{{a + d}}{{u + d}}} \right)} \right)}_ + },\log b} \right\}.\]
Then, for any ${B_r} \equiv {B_r}\left( {{\xi _0}} \right) \subset {B_{\frac{R}{2}}}\left( {{\xi _0}} \right),\;r \in \left( {0,1} \right]$, there exists a positive constant $c = c\left( {n,p,s} \right)$ such that
\begin{equation}\label{eq46}
\int_{{B_r}} {{{\left| {v - {{\left( v \right)}_r}} \right|}^p}{\mkern 1mu} d\xi }  \le c\left( {1 + {d^{1 - p}}{{\left( {\frac{r}{R}} \right)}^p}{\rm{Tail}}{{({u_ - };{\xi _0},R)}^{p - 1}}} \right),
\end{equation}
where ${\rm{Tail}}\left(  \cdot  \right)$ is given by \eqref{eq21}.
\end{corollary}

\begin{lemma}\label{Le43}
Let $u \in HW_{loc}^{1,p}\left( \Omega  \right)\;\;(1 < p < Q)$ be a weak solution to \eqref{eq15}, and ${B_R}\left( {{\xi _0}} \right) \subset \Omega ,\;r \in \left( {0,1} \right],\;r \in \left( {0,\frac{R}{2}} \right)$. For any $\sigma  \in \left( {0,\frac{1}{4}} \right]$, denote
\[{r_j} = {\sigma ^j}\frac{r}{2},\;{B_j} \equiv {B_{{r_j}}}\left( {{\xi _0}} \right),\;j = 0,1,2, \cdots ,\]
\begin{equation}\label{eq47}
\frac{1}{2}\vartheta \left( {{r_0}} \right): = \frac{1}{2}\vartheta \left( {\frac{r}{2}} \right) = {\rm{Tail}}(u;{\xi _0},\frac{r}{2}) + c{\left( {\fint_{{B_r}} {{{\left| u \right|}^p}d\xi } } \right)^{\frac{1}{p}}},
\end{equation}
where ${\rm{Tail}}\left(  \cdot  \right)$ is given by \eqref{eq21}, and $c = c\left( {n,p,s} \right)$ is the constant in \eqref{eq16}. For $\alpha  \in \left( {0\frac{p}{{p - 1}}} \right)$, if we let
\begin{equation}\label{eq48}
\vartheta \left( {{r_j}} \right) = {\left( {\frac{{{r_j}}}{{{r_0}}}} \right)^\alpha }\vartheta \left( {{r_0}} \right),\;j = 0,1,2, \cdots ,
\end{equation}
then
\begin{equation}\label{eq49}
\mathop {{\rm{osc }}u}\limits_{{B_j}}  = \mathop {{\rm{ess}}\sup }\limits_{{B_j}} u - \mathop {{\rm{ess inf}}}\limits_{{B_j}} u \le \vartheta \left( {{r_j}} \right),\;j = 0,1,2, \cdots .
\end{equation}
\end{lemma}

\begin{proof}
We proceed with the proof by the induction.

\textbf{Step 1.} Since $u$ is a weak solution to \eqref{eq15}, both $u$ and $-u$ are weak subsolutions to \eqref{eq15}. Thus by Lemma \ref{Le25}, ${u_ + }$ and ${\left( { - u} \right)_ + }$ are weak subsolutions to \eqref{eq15}. Therefore, applying Theorem \ref{Th11} to ${u_ + }$ and ${\left( { - u} \right)_ + }$, and taking $\delta  = 1$, we obtain that \eqref{eq49} holds for $j = 0$.

\textbf{Step 2.} Suppose that \eqref{eq49} holds for some $j \in \left\{ {1,2, \cdots } \right\}$; we next prove that \eqref{eq49} holds for $j + 1$. Note that either
\begin{equation}\label{eq410}
\frac{{\left| {{B_{2{r_{j + 1}}}}\left( {{\xi _0}} \right) \cap \left\{ {u \ge \mathop {{\rm{ess inf}}}\limits_{{B_j}} u + \frac{{\vartheta \left( {{r_j}} \right)}}{2}} \right\}} \right|}}{{{B_{2{r_{j + 1}}}}\left( {{\xi _0}} \right)}} \ge \frac{1}{2}
\end{equation}
or
\begin{equation}\label{eq411}
\frac{{\left| {{B_{2{r_{j + 1}}}}\left( {{\xi _0}} \right) \cap \left\{ {u \le \mathop {{\rm{ess inf}}}\limits_{{B_j}} u + \frac{{\vartheta \left( {{r_j}} \right)}}{2}} \right\}} \right|}}{{{B_{2{r_{j + 1}}}}\left( {{\xi _0}} \right)}} \ge \frac{1}{2}
\end{equation}
must hold. Let
\begin{equation}\label{eq412}
{u_j} = \left\{ \begin{array}{l}
u - \mathop {{\rm{ess inf}}}\limits_{{B_j}} u,\;\;\;\;\;\;\;\;\;\;\;\;\;\;\;\;\;{\rm{if}} \;\eqref{eq410} \;{\rm{holds}},\\
\vartheta \left( {{r_j}} \right) - \left( {u - \mathop {{\rm{ess inf}}}\limits_{{B_j}} u} \right),\;\;{\mathop{\rm if}\nolimits} \; \eqref{eq411}\; {\rm{holds}}.
\end{array} \right.
\end{equation}
Then ${u_j}$ is a weak solution to \eqref{eq15}, and whether \eqref{eq410} or \eqref{eq411} holds, we have ${u_j} \ge 0$ in ${B_j}$ and
\begin{equation}\label{eq40}
\frac{{\left| {{B_{2{r_{j + 1}}}}\left( {{\xi _0}} \right) \cap \left\{ {{u_j} \ge \frac{{\vartheta \left( {{r_j}} \right)}}{2}} \right\}} \right|}}{{{B_{2{r_{j + 1}}}}\left( {{\xi _0}} \right)}} \ge \frac{1}{2}.
\end{equation}

Denote
\begin{equation}\label{eq413}
v = \min \left\{ {{{\left( {\log \left( {\frac{{\frac{{\vartheta \left( {{r_j}} \right)}}{2} + d}}{{{u_j} + d}}} \right)} \right)}_ + },k} \right\},\;\;k > 0.
\end{equation}

\textbf{Claim 1:} It holds
\begin{equation}\label{eq414}
\fint_{{B_{2{r_{j + 1}}}}} {{{\left| {v - {{\left( v \right)}_{{B_{2{r_{j + 1}}}}}}} \right|}^p}{\mkern 1mu} d\xi }  \le C.
\end{equation}

In fact, since $\sigma  \in \left( {0,\frac{1}{4}} \right]$, we have $2{r_{j + 1}} = 2{\sigma ^{j + 1}}\frac{r}{2} \le \frac{1}{2}{\sigma ^j}\frac{r}{2} = \frac{{{r_j}}}{2}$ and thus ${B_{2{r_{j + 1}}}} \subset {B_{\frac{{{r_j}}}{2}}} \subset {B_{{r_j}}}$. Applying Corollary \ref{Co42} to $v$ and taking $a = \frac{{\vartheta \left( {{r_j}} \right)}}{2}$ and $b = {e^k}$, we obtain
\begin{align}\label{eq415}
   \fint_{{B_{2{r_{j + 1}}}}} {{{\left| {v - {{\left( v \right)}_{{B_{2{r_{j + 1}}}}}}} \right|}^p}{\mkern 1mu} d\xi }& \le c\left( {1 + {d^{1 - p}}{{\left( {\frac{{2{r_{j + 1}}}}{{{r_j}}}} \right)}^p}{\rm{Tail}}{{({u_j};{\xi _0},{r_j})}^{p - 1}}} \right)\nonumber \\
   &  \le c\left( {1 + {d^{1 - p}}{\sigma ^p}{\rm{Tail}}{{({u_j};{\xi _0},{r_j})}^{p - 1}}} \right).
\end{align}
Moreover, it follows from (6.6) in [\cite{MPPP23}, p32] that
\begin{equation}\label{eq416}
{\left[ {{\rm{Tail}}({u_j};{\xi _0},{r_j})} \right]^{p - 1}} \le c\left( {n,p,s,\alpha } \right){\sigma ^{ - \alpha \left( {p - 1} \right)}}{\left[ {\vartheta \left( {{r_j}} \right)} \right]^{p - 1}}.
\end{equation}
Taking $d = {\sigma ^{\frac{p}{{p - 1}}}}$ and substituting \eqref{eq416} into \eqref{eq415}, it yields
\begin{align}\label{eq417}
   \fint_{{B_{2{r_{j + 1}}}}} {{{\left| {v - {{\left( v \right)}_{{B_{2{r_{j + 1}}}}}}} \right|}^p}{\mkern 1mu} d\xi }& \le c\left( {1 + {d^{1 - p}}{\sigma ^p}{\sigma ^{ - \alpha \left( {p - 1} \right)}}{{\left[ {\vartheta \left( {{r_j}} \right)} \right]}^{p - 1}}} \right)\nonumber \\
   &  = c\left( {1 + {\sigma ^{\alpha \left( {p - 1} \right)\left( {j - 1} \right)}}{{\left[ {\vartheta \left( {{r_0}} \right)} \right]}^{p - 1}}} \right) \le C.
\end{align}

\textbf{Claim 2:} It holds
\begin{equation}\label{eq418}
\frac{{\left| {{B_{2{r_{j + 1}}}} \cap \left\{ {{u_j} \le 2d\vartheta \left( {{r_j}} \right)} \right\}} \right|}}{{\left| {{B_{2{r_{j + 1}}}}} \right|}} \le \frac{{\tilde c\left( {n,p,s} \right)}}{{\log \left( {\frac{1}{\sigma }} \right)}}.
\end{equation}

Actually, similar to the proof of (6.13) in [\cite{MPPP23}, p35], by using \eqref{eq40}-\eqref{eq414}, we can gain \eqref{eq418}.

Next, similar to the proof of Theorem \ref{Th11}, we can prove that \eqref{eq49} holds for $j+1$ by iteration.

In fact, for any $i = 0,1,2, \cdots $, define
\[{\rho _i} = {r_{j + 1}} + {2^{ - i}}{r_{j + 1}},\;{\bar \rho _i} = \frac{{{\rho _{i + 1}} + {\rho _i}}}{2},\;{B^i} = {B_{{\rho _i}}},\;{\bar B^i} = {B_{{{\bar \rho }_i}}};\]
\[{\psi _i} \in C_0^\infty \left( {{{\bar B}^i}} \right),\;0 \le {\psi _i} \le 1,\;\left| {{\nabla _H}{\psi _i}} \right| \le \frac{c}{{{\rho _i}}}\;\rm{in} \; {{{\bar B}^i}}, {\psi _i} = 1\;{\rm{in}}\;{B^{i + 1}};\]
\[{k_i} = \left( {1 + {2^{ - i}}} \right)d\vartheta \left( {{r_j}} \right),\;{\omega _i} = {\left( {{k_i} - {u_j}} \right)_ + };\]
\[{A_i} = \frac{{\left| {{B^i} \cap \left\{ {{u_j} \le {k_i}} \right\}} \right|}}{{\left| {{B^i}} \right|}} = \frac{{\left| {{B^i} \cap \left\{ {{\omega _i} > 0} \right\}} \right|}}{{\left| {{B^i}} \right|}}.\]
It follows from \eqref{eq418} that
\begin{equation}\label{eq419}
{A_0} = \frac{{\left| {{B_{2{r_{j + 1}}}} \cap \left\{ {{u_j} \le 2d\vartheta \left( {{r_j}} \right)} \right\}} \right|}}{{\left| {{B_{2{r_{j + 1}}}}} \right|}} \le \frac{{\tilde c\left( {n,p,s} \right)}}{{\log \left( {\frac{1}{\sigma }} \right)}}.
\end{equation}

On the one hand, by \eqref{eq23}, we have
\begin{align}\label{eq420}
   & \int_{{B^i}} {\psi _i^p{{\left| {{\nabla _H}{\omega _i}} \right|}^p}d\xi } \nonumber \\
  \le &  \int_{{B^i}} {\psi _i^p{{\left| {{\nabla _H}{\omega _i}} \right|}^p}d\xi }  + \int_{{B^i}} {\int_{{B^i}} {\frac{{{{\left| {{\omega _i}\left( \xi  \right){\psi _i}\left( \xi  \right) - {\omega _i}\left( \eta  \right){\psi _i}\left( \eta  \right)} \right|}^p}}}{{\left\| {{\eta ^{ - 1}} \circ \xi } \right\|_{{{\rm{H}}^n}}^{Q + sp}}}{\mkern 1mu} d\xi d\eta } } \nonumber \\
    \le& c(\int_{{B^i}} {\omega _i^p{{\left| {{\nabla _H}{\psi _i}} \right|}^p}d\xi }  + \int_{{B^i}} {\int_{{B^i}} {\max {{\{ {\omega _i}\left( \xi  \right),{\omega _i}\left( \eta  \right)\} }^p}\frac{{{{\left| {{\psi _i}\left( \xi  \right) - {\psi _i}\left( \eta  \right)} \right|}^p}}}{{\left\| {{\eta ^{ - 1}} \circ \xi } \right\|_{{{\rm{H}}^n}}^{Q + sp}}}{\mkern 1mu} d\xi d\eta } } \nonumber \\
  & + \mathop {{\rm{ess}}\;\sup }\limits_{\xi  \in {\rm{supp}}\;{\psi _i}} \int_{{{\rm{H}}^n}\backslash {B^i}} {\frac{{{\omega _i}{{\left( \eta  \right)}^{p - 1}}}}{{\left\| {{\eta ^{ - 1}} \circ \xi } \right\|_{{{\rm{H}}^n}}^{Q + sp}}}d\eta }  \cdot \int_{{B^i}} {{\omega _i}\psi _i^pd\xi } ).
\end{align}
On the other hand, it deduces from Lemma \ref{Le22} that
\begin{align}\label{eq421}
   A_{i + 1}^{\frac{p}{{{p^*}}}}{\left( {{k_i} - {k_{i + 1}}} \right)^p}&= \frac{{{{\left| {{B^{i + 1}} \cap \left\{ {{u_j} \le {k_{i + 1}}} \right\}} \right|}^{\frac{p}{{{p^*}}}}}}}{{{{\left| {{B^{i + 1}}} \right|}^{\frac{p}{{{p^*}}}}}}}{\left( {{k_i} - {k_{i + 1}}} \right)^p} \nonumber \\
   & = \frac{1}{{{{\left| {{B^{i + 1}}} \right|}^{\frac{p}{{{p^*}}}}}}}{\left( {\int_{{B^{i + 1}} \cap \left\{ {{u_j} \le {k_{i + 1}}} \right\}} {{{\left( {{k_i} - {k_{i + 1}}} \right)}^{{p^*}}}\psi _i^{{p^*}}d\xi } } \right)^{\frac{p}{{{p^*}}}}} \nonumber \\
   &\le \frac{{{{\left| {{B^i}} \right|}^{\frac{p}{{{p^*}}}}}}}{{{{\left| {{B^{i + 1}}} \right|}^{\frac{p}{{{p^*}}}}}}}{\left( {\fint_{{B^i}} {\omega _i^{{p^*}}\psi _i^{{p^*}}d\xi } } \right)^{\frac{p}{{{p^*}}}}}\nonumber \\
  & \le cr_{j + 1}^p\fint_{{B^i}} {{{\left| {{\nabla _H}\left( {\omega _i\psi _i} \right)} \right|}^p}d\xi }  \nonumber \\
  & = cr_{j + 1}^p\fint_{{B^i}} {\psi _i^p{{\left| {{\nabla _H}\omega _i} \right|}^p}d\xi }  + cr_{j + 1}^p\fint_{{B^i}} {\omega _i^p{{\left| {{\nabla _H}\psi _i} \right|}^p}d\xi } ,
\end{align}
and so combining \eqref{eq420} and \eqref{eq421}, we derive
\begin{align}\label{eq422}
   & A_{i + 1}^{\frac{p}{{{p^*}}}}{\left( {{k_i} - {k_{i + 1}}} \right)^p}\nonumber \\
   \le& cr_{j + 1}^p(\fint_{{B^i}} {\omega _i^p{{\left| {{\nabla _H}{\psi _i}} \right|}^p}d\xi }  + \int_{{B^i}} {\fint_{{B^i}} {\max {{\{ {\omega _i}\left( \xi  \right),{\omega _i}\left( \eta  \right)\} }^p}\frac{{{{\left| {{\psi _i}\left( \xi  \right) - {\psi _i}\left( \eta  \right)} \right|}^p}}}{{\left\| {{\eta ^{ - 1}} \circ \xi } \right\|_{{{\rm{H}}^n}}^{Q + sp}}}{\mkern 1mu} d\xi d\eta } }\nonumber \\
    &+ \mathop {{\rm{ess}}\;\sup }\limits_{\xi  \in {\rm{supp}}\;{\psi _i}} \int_{{{\rm{H}}^n}\backslash {B^i}} {\frac{{{\omega _i}{{\left( \eta  \right)}^{p - 1}}}}{{\left\| {{\eta ^{ - 1}} \circ \xi } \right\|_{{{\rm{H}}^n}}^{Q + sp}}}d\eta }  \cdot\fint_{{B^i}} {{\omega _i}\psi _i^pd\xi } )\nonumber \\
     = :&cr_{j + 1}^p\left( {{I_1} + {I_2} + {I_3}} \right).
\end{align}

By using $\left| {{\nabla _H}{\psi _i}} \right| \le \frac{c}{{{\rho _i}}}$ and ${\omega _i} \le {k_i} \le 2d\vartheta \left( {{r_j}} \right)$, it implies
\begin{equation}\label{eq423}
{I_1} = \fint_{{B^i}} {\omega _i^p{{\left| {{\nabla _H}{\psi _i}} \right|}^p}d\xi }  \le c{2^{ip}}r_{j + 1}^{ - p}{\left( {d\vartheta \left( {{r_j}} \right)} \right)^p}{A_i}.
\end{equation}
It yields from (6.16) in [\cite{MPPP23}, p36] that
\begin{equation}\label{eq424}
{I_2} \le c{2^{ip}}r_{j + 1}^{ - sp}{\left( {d\vartheta \left( {{r_j}} \right)} \right)^p}{A_i} \le c{2^{ip}}r_{j + 1}^{ - p}{\left( {d\vartheta \left( {{r_j}} \right)} \right)^p}{A_i}.
\end{equation}
It follows from (6.17) and (6.19) in [\cite{MPPP23}, pp36-37] that
\begin{equation}\label{eq425}
{I_3} \le c{2^{i\left( {Q + sp} \right)}}r_{j + 1}^{ - sp}{\left( {d\vartheta \left( {{r_j}} \right)} \right)^p}{A_i} \le c{2^{i\left( {Q + p} \right)}}r_{j + 1}^{ - p}{\left( {d\vartheta \left( {{r_j}} \right)} \right)^p}{A_i}.
\end{equation}
Substituting \eqref{eq423}-\eqref{eq425} into \eqref{eq422}, we have
\[A_{i + 1}^{\frac{p}{{{p^*}}}}{\left( {{k_i} - {k_{i + 1}}} \right)^p} \le c{2^{i\left( {Q + p} \right)}}{\left( {d\vartheta \left( {{r_j}} \right)} \right)^p}{A_i},\]
so
\[{A_{i + 1}} \le c{2^{i\left( {Q + 2p} \right)\frac{{{p^*}}}{p}}}A_i^{\frac{{{p^*}}}{p}} = c{2^{i\left( {Q + 2p} \right)\frac{{{p^*}}}{p}}}A_i^{1 + \frac{p}{{Q - p}}}.\]
If
\begin{equation}\label{eq426}
{A_0} \le {c^{ - \frac{{Q - p}}{p}}}{2^{ - \left( {Q + 2p} \right)\frac{{{p^*}}}{p}{{\left( {\frac{{Q - p}}{p}} \right)}^2}}} \le {c^{ - \frac{{Q - p}}{p}}}{2^{ - \left( {Q + 2p} \right)\frac{{Q\left( {Q - p} \right)}}{{{p^2}}}}} = :v*,
\end{equation}
then by Lemma \ref{Le26}, we gain
\[\mathop {\lim }\limits_{i \to \infty } {A_i} = \mathop {\lim }\limits_{i \to \infty } \frac{{\left| {{B^i} \cap \left\{ {{\omega _i} > 0} \right\}} \right|}}{{\left| {{B^i}} \right|}} = 0,\]
i.e.
\[\mathop {{\rm{ess}}\;{\rm{inf}}}\limits_{{B_{j + 1}}} \;{u_j} \ge d\vartheta \left( {{r_j}} \right).\]
Thus, by Theorem \ref{Th11}, \eqref{eq412} and the above formula, we have
\[\mathop {{\rm{osc}}u}\limits_{{B_{j + 1}}}  = \mathop {{\rm{ess}}\sup }\limits_{{B_{j + 1}}} {u_j} - \mathop {{\rm{ess inf}}}\limits_{{B_{j + 1}}} {u_j} \le \left( {1 - d} \right)\vartheta \left( {{r_j}} \right) = \left( {1 - d} \right){\left( {\frac{{{r_j}}}{{{r_{j + 1}}}}} \right)^\alpha }\vartheta \left( {{r_{j + 1}}} \right) = \left( {1 - d} \right){\sigma ^{ - \alpha }}\vartheta \left( {{r_{j + 1}}} \right).\]
If we take $\alpha  \in \left( {0,\frac{p}{{p - 1}}} \right)$ sufficiently small such that
\[{\sigma ^\alpha } \ge 1 - d = 1 - {\sigma ^{\frac{p}{{p - 1}}}},\]
then
\[\mathop {{\rm{osc}}u}\limits_{{B_{j + 1}}}  \le \vartheta \left( {{r_{j + 1}}} \right).\]
In fact, if we take $\sigma  = \min \left\{ {\frac{1}{4},{e^{ - \frac{{\tilde c\left( {n,p,s} \right)}}{{v*}}}}} \right\}$ and combine it with \eqref{eq419}, we can obtain that \eqref{eq426} holds. Therefore, this lemma is established.
\end{proof}

\textbf{Proof of Theorem \ref{Th12}.} Theorem \ref{Th12} can be directly obtained from Lemma \ref{Le43}.

\section{Harnack Inequality}
\label{Section 5}

The main purpose of this section is to prove Theorem \ref{Th13} by using the iterative lemma (Lemma \ref{Le54}) and an integral estimate (Lemma \ref{Le53}). To apply Lemma \ref{Le54}, it is necessary to first prove that the condition \eqref{eq531} holds. Therefore, we first present a Tail estimate (Lemma \ref{Le51}), and then prove that \eqref{eq531} holds by using Theorem \ref{Th11} and Lemma \ref{Le51}. On the other hand, to prove Lemma \ref{Le53}, we inspect a expansion of positivity (Lemma \ref{Le52}). In addition, Lemma \ref{Le53} will also be used in the next section to prove the weak Harnack inequality.

\begin{lemma}[Tail estimate]\label{Le51}
Let $u \in HW_{loc}^{1,p}\left( \Omega  \right)\;\left( {1 < p < \infty } \right)$ be a weak solution to \eqref{eq15} and $u \ge 0$ in ${B_R}\left( {{\xi _0}} \right) \subset \Omega $. Then, for $0 < r < R,\;r \in \left( {0,1} \right]$, there exists a positive constant $c = c\left( {n,p,s} \right)$ such that
\begin{equation}\label{eq51}
{\rm{Tail}}({u_ + };{\xi _0},r) \le c\mathop {{\rm{ess}}\sup }\limits_{{B_r}\left( {{\xi _0}} \right)} u + c{\left( {\frac{r}{R}} \right)^{\frac{p}{{p - 1}}}}{\rm{Tail}}({u_ - };{\xi _0},R),
\end{equation}
where ${\rm{Tail}}\left(  \cdot  \right)$ is given by \eqref{eq21}.
\end{lemma}

\begin{proof}
Let $M = \mathop {{\rm{ess}}\sup }\limits_{{B_r}\left( {{\xi _0}} \right)} u$ and $\psi  \in C_0^\infty \left( {{B_r}\left( {{\xi _0}} \right)} \right)$ be a cut-off function satisfying
\begin{center}
$0 \le \psi  \le 1,\;\left| {{\nabla _H}\psi } \right| \le \frac{8}{r}$ in ${B_r}\left( {{\xi _0}} \right)$ and $\psi  = 1$ in ${B_{\frac{r}{2}}}\left( {{\xi _0}} \right)$.
\end{center}
Denote $\omega  = u - 2M$ and take $\phi  = \omega {\psi ^p}$ as the test function in \eqref{eq22}. Then
\begin{align}\label{eq52}
  0 = & \int_{{B_r}} {{{\left| {{\nabla _H}u} \right|}^{p - 2}}{\nabla _H}u \cdot {\nabla _H}\left( {\omega {\psi ^p}} \right)d\xi }  \nonumber\\
   &  + \int_{{{\rm{H}}^n}} {\int_{{{\rm{H}}^n}} {\frac{{{{\left| {u\left( \xi  \right) - u\left( \eta  \right)} \right|}^{p - 2}}\left( {u\left( \xi  \right) - u\left( \eta  \right)} \right)\left( {\omega \left( \xi  \right)\psi {{\left( \xi  \right)}^p} - \omega \left( \eta  \right)\psi {{\left( \eta  \right)}^p}} \right)}}{{\left\| {{\eta ^{ - 1}} \circ \xi } \right\|_{{{\rm{H}}^n}}^{Q + sp}}}d\xi d\eta } } \nonumber\\
   = &\int_{{B_r}} {{{\left| {{\nabla _H}\omega } \right|}^{p - 2}}{\nabla _H}\omega  \cdot {\nabla _H}\left( {\omega {\psi ^p}} \right)d\xi } \nonumber\\
  & + \int_{{B_r}} {\int_{{B_r}} {\frac{{{{\left| {\omega \left( \xi  \right) - \omega \left( \eta  \right)} \right|}^{p - 2}}\left( {\omega \left( \xi  \right) - \omega \left( \eta  \right)} \right)\left( {\omega \left( \xi  \right)\psi {{\left( \xi  \right)}^p} - \omega \left( \eta  \right)\psi {{\left( \eta  \right)}^p}} \right)}}{{\left\| {{\eta ^{ - 1}} \circ \xi } \right\|_{{{\rm{H}}^n}}^{Q + sp}}}d\xi d\eta } } \nonumber\\
   & + 2\int_{{{\rm{H}}^n}\backslash {B_r}} {\int_{{B_r}} {\frac{{{{\left| {\omega \left( \xi  \right) - \omega \left( \eta  \right)} \right|}^{p - 2}}\left( {\omega \left( \xi  \right) - \omega \left( \eta  \right)} \right)\omega \left( \xi  \right)\psi {{\left( \xi  \right)}^p}}}{{\left\| {{\eta ^{ - 1}} \circ \xi } \right\|_{{{\rm{H}}^n}}^{Q + sp}}}d\xi d\eta } } \nonumber\\
   =:& {I_1} + {I_2} +{I_3}.
\end{align}

For ${I_1}$, by using Young's inequality, we obtain
\begin{align}\label{eq53}
   {I_1}& = \int_{{B_r}} {{{\left| {{\nabla _H}\omega } \right|}^{p - 2}}{\nabla _H}\omega  \cdot {\nabla _H}\left( {\omega {\psi ^p}} \right)d\xi }  \nonumber\\
   &  = \int_{{B_r}} {{{\left| {{\nabla _H}\omega } \right|}^p}{\psi ^p}d\xi }  + p\int_{{B_r}} {{\psi ^{p - 1}}\omega {{\left| {{\nabla _H}\omega } \right|}^{p - 2}}{\nabla _H}\omega  \cdot {\nabla _H}\psi d\xi } \nonumber\\
   & \ge \frac{1}{2}\int_{{B_r}} {{{\left| {{\nabla _H}\omega } \right|}^p}{\psi ^p}d\xi }  - c\left( p \right)\int_{{B_r}} {{{\left| \omega  \right|}^p}{{\left| {{\nabla _H}\psi } \right|}^p}d\xi } \nonumber\\
  & \ge  - c\left( p \right)\int_{{B_r}} {{{\left| M \right|}^p}{{\left| {{\nabla _H}\psi } \right|}^p}d\xi }  \nonumber\\
   & \ge  - c\left( p \right){\left| M \right|^p}{r^{ - p}}\left| {{B_r}} \right|.
\end{align}
For ${I_2}$, it follows from the derivation of (4.24) in \cite{PP22} that
\begin{equation}\label{eq54}
{I_2} \ge  - c{M^p}{r^{ - sp}}\left| {{B_r}} \right| \ge  - c{M^p}{r^{ - p}}\left| {{B_r}} \right|.
\end{equation}
For ${I_3}$, it deduces from the derivation of (4.23) in \cite{PP22} that
\begin{align}\label{eq55}
   {I_3}& \ge cM{r^{ - p}}\left| {{B_r}} \right|{\left[ {{\rm{Tail}}({u_ + };{\xi _0},r)} \right]^{p - 1}} - c{M^p}{r^{ - sp}}\left| {{B_r}} \right| - cM\left| {{B_r}} \right|{R^{ - p}}{\left[ {{\rm{Tail}}({u_ - };{\xi _0},R)} \right]^{p - 1}}\nonumber \\
   &  \ge cM{r^{ - p}}\left| {{B_r}} \right|{\left[ {{\rm{Tail}}({u_ + };{\xi _0},r)} \right]^{p - 1}} - c{M^p}{r^{ - p}}\left| {{B_r}} \right| - cM\left| {{B_r}} \right|{R^{ - p}}{\left[ {{\rm{Tail}}({u_ - };{\xi _0},R)} \right]^{p - 1}}.
\end{align}

Therefore, substituting \eqref{eq53}-\eqref{eq55} into \eqref{eq52} yields \eqref{eq51}.
\end{proof}

\begin{lemma}[Expansion of positivity]\label{Le52}
Let $u \in HW_{loc}^{1,p}\left( \Omega  \right)\;\left( {1 < p < Q } \right)$ be a weak supersolution to \eqref{eq15} and $u \ge 0$ in ${B_R}\left( {{\xi _0}} \right) \subset \Omega $. If for $k \ge 0,$ $r \in \left( {0,1} \right],$ $0 < r < \frac{R}{{16}}$, there exists $\tau  \in \left( {0,1} \right]$ such that
\begin{equation}\label{eq56}
\left| {{B_r}\left( {{\xi _0}} \right) \cap \left\{ {u \ge k} \right\}} \right| \ge \tau \left| {{B_r}\left( {{\xi _0}} \right)} \right|,
\end{equation}
then there exists a constant $\delta  = \delta \left( {n,p,s,\tau } \right) \in \left( {0,\frac{1}{4}} \right)$ such that
\begin{equation}\label{eq57}
\mathop {{\rm{ess inf}}}\limits_{{B_{4r}}\left( {{\xi _0}} \right)} u \ge \delta k - {\left( {\frac{r}{R}} \right)^{\frac{p}{{p - 1}}}}{\rm{Tail}}({u_ - };{\xi _0},R),
\end{equation}
where ${\rm{Tail}}\left(  \cdot  \right)$ is given by \eqref{eq21}.
\end{lemma}

\begin{proof}
The proof is completed by two steps.

\textbf{Step 1.} Under the assumption \eqref{eq56}, we prove that there exists a constant ${c_1} = c\left( {n,p,s} \right)$ such that for every $\delta  \in \left( {0,\frac{1}{4}} \right)$ and $\varepsilon  > 0$, the following holds:
\begin{equation}\label{eq58}
\left| {{B_{6r}}\left( {{\xi _0}} \right) \cap \left\{ {u \le 2\delta k - \frac{1}{2}{{\left( {\frac{r}{R}} \right)}^{\frac{p}{{p - 1}}}}{\rm{Tail}}({u_ - };{\xi _0},R) - \varepsilon } \right\}} \right| \le \frac{{{c_1}}}{{\tau \log \frac{1}{{2\delta }}}}\left| {{B_{6r}}\left( {{\xi _0}} \right)} \right|.
\end{equation}

In fact, let $\varepsilon  > 0$ and $\psi  \in C_0^\infty \left( {{B_{7r}}\left( {{\xi _0}} \right)} \right)$ be a cut-off function satisfying
\begin{center}
$0 \le \psi  \le 1,\;\left| {{\nabla _H}\psi } \right| \le \frac{8}{r}$ in ${B_{7r}}\left( {{\xi _0}} \right)$ and $\psi  = 1$ in ${B_{6r}}\left( {{\xi _0}} \right)$.
\end{center}
Denote
\[\omega  = u + {d_\varepsilon },\]
where
\[{d_\varepsilon } = \frac{1}{2}{\left( {\frac{r}{R}} \right)^{\frac{p}{{p - 1}}}}{\rm{Tail}}({u_ - };{\xi _0},R) + \varepsilon .\]
Taking $\phi  = {\omega ^{1 - p}}{\psi ^p}$ as the test function in \eqref{eq22}, we obtain
\begin{align}\label{eq59}
   0 \le& \int_{{B_{8r}}} {{{\left| {{\nabla _H}u} \right|}^{p - 2}}{\nabla _H}u \cdot {\nabla _H}\left( {{\omega ^{1 - p}}{\psi ^p}} \right)d\xi } \nonumber \\
   &  + \int_{{{\rm{H}}^n}} {\int_{{{\rm{H}}^n}} {\frac{{{{\left| {u\left( \xi  \right) - u\left( \eta  \right)} \right|}^{p - 2}}\left( {u\left( \xi  \right) - u\left( \eta  \right)} \right)\left( {\omega {{\left( \xi  \right)}^{1 - p}}\psi {{\left( \xi  \right)}^p} - \omega {{\left( \eta  \right)}^{1 - p}}\psi {{\left( \eta  \right)}^p}} \right)}}{{\left\| {{\eta ^{ - 1}} \circ \xi } \right\|_{{{\rm{H}}^n}}^{Q + sp}}}d\xi d\eta } } \nonumber \\
    =& \int_{{B_{8r}}} {{{\left| {{\nabla _H}\omega } \right|}^{p - 2}}{\nabla _H}\omega  \cdot {\nabla _H}\left( {{\omega ^{1 - p}}{\psi ^p}} \right)d\xi } \nonumber \\
   & + \int_{{B_{8r}}} {\int_{{B_{8r}}} {\frac{{{{\left| {\omega \left( \xi  \right) - \omega \left( \eta  \right)} \right|}^{p - 2}}\left( {\omega \left( \xi  \right) - \omega \left( \eta  \right)} \right)\left( {\omega {{\left( \xi  \right)}^{1 - p}}\psi {{\left( \xi  \right)}^p} - \omega {{\left( \eta  \right)}^{1 - p}}\psi {{\left( \eta  \right)}^p}} \right)}}{{\left\| {{\eta ^{ - 1}} \circ \xi } \right\|_{{{\rm{H}}^n}}^{Q + sp}}}d\xi d\eta } } \nonumber \\
   & + 2\int_{{{\rm{H}}^n}\backslash {B_{8r}}} {\int_{{B_{8r}}} {\frac{{{{\left| {\omega \left( \xi  \right) - \omega \left( \eta  \right)} \right|}^{p - 2}}\left( {\omega \left( \xi  \right) - \omega \left( \eta  \right)} \right)\omega {{\left( \xi  \right)}^{1 - p}}\psi {{\left( \xi  \right)}^p}}}{{\left\| {{\eta ^{ - 1}} \circ \xi } \right\|_{{{\rm{H}}^n}}^{Q + sp}}}d\xi d\eta } } \nonumber \\
    =:& {I_1} + {I_2} +{I_3}.
\end{align}

For ${I_1}$, similar to the estimate of \eqref{eq43}, we have
\begin{align}\label{eq510}
   {I_1}& = \int_{{B_{8r}}} {{{\left| {{\nabla _H}\omega } \right|}^{p - 2}}{\nabla _H}\omega  \cdot {\nabla _H}\left( {{\omega ^{1 - p}}{\psi ^p}} \right)d\xi } \nonumber \\
   &  \le  - c\left( p \right)\int_{{B_{8r}}} {{{\left| {{\nabla _H}\log \omega } \right|}^p}{\psi ^p}d\xi }  + c\left( p \right){r^{Q - p}}\nonumber \\
   & \le  - c\left( p \right)\int_{{B_{6r}}} {{{\left| {{\nabla _H}\log \omega } \right|}^p}{\psi ^p}d\xi }  + c\left( p \right){r^{Q - p}}\nonumber \\
   & \le  - c\left( p \right)\int_{{B_{6r}}} {{{\left| {{\nabla _H}\log \omega } \right|}^p}d\xi }  + c\left( p \right){r^{Q - p}}.
\end{align}
For ${I_2}$, it follows from the estimate of ${I_1}$ in Lemma 4.1 of \cite{PP22} that
\begin{align}\label{eq511}
  {I_2} & \le  - \frac{1}{c}\int_{{B_{6r}}} {\int_{{B_{6r}}} {{{\left| {\log \left( {\frac{{\omega \left( \xi  \right)}}{{\omega \left( \eta  \right)}}} \right)} \right|}^p}\frac{1}{{\left\| {{\eta ^{ - 1}} \circ \xi } \right\|_{{{\rm{H}}^n}}^{Q + sp}}}d\xi d\eta } }  + c{r^{Q - sp}}\nonumber \\
  &  \le  - \frac{1}{c}\int_{{B_{6r}}} {\int_{{B_{6r}}} {{{\left| {\log \left( {\frac{{\omega \left( \xi  \right)}}{{\omega \left( \eta  \right)}}} \right)} \right|}^p}\frac{1}{{\left\| {{\eta ^{ - 1}} \circ \xi } \right\|_{{{\rm{H}}^n}}^{Q + sp}}}d\xi d\eta } }  + c{r^{Q - p}}.
\end{align}
For ${I_3}$, it holds from the estimates of ${I_2}$ and ${I_3}$ in Lemma 4.1 of \cite{PP22} that
\begin{equation}\label{eq512}
{I_3} \le c{r^{Q - sp}} \le c{r^{Q - p}}.
\end{equation}

Therefore, substituting \eqref{eq510}-\eqref{eq512} into \eqref{eq59} yields
\begin{equation}\label{eq513}
\int_{{B_{6r}}} {{{\left| {{\nabla _H}\log \omega } \right|}^p}d\xi }  + \int_{{B_{6r}}} {\int_{{B_{6r}}} {{{\left| {\log \left( {\frac{{\omega \left( \xi  \right)}}{{\omega \left( \eta  \right)}}} \right)} \right|}^p}\frac{1}{{\left\| {{\eta ^{ - 1}} \circ \xi } \right\|_{{{\rm{H}}^n}}^{Q + sp}}}d\xi d\eta } }  \le c{r^{Q - p}}.
\end{equation}

For $\delta  \in \left( {0,\frac{1}{4}} \right)$, denote
\[\upsilon  = {\left( {\min \left\{ {\log \frac{1}{{2\delta }},\log \frac{{k + {d_\varepsilon }}}{\omega }} \right\}} \right)_ + }.\]
Then, by H\"{o}lder's inequality, Lemma \ref{Le21} and \eqref{eq513}, we have
\begin{align}\label{eq514}
   \int_{{B_{6r}}} {\left| {\upsilon \left( \xi  \right) - {{\left( \upsilon  \right)}_{{B_{6r}}}}} \right|d\xi }& \le {\left| {{B_{6r}}} \right|^{1 - \frac{{Q - p}}{{pQ}}}}{\left( {\int_{{B_{6r}}} {{{\left| {\upsilon \left( \xi  \right) - {{\left( \upsilon  \right)}_{{B_{6r}}}}} \right|}^{\frac{{pQ}}{{Q - p}}}}d\xi } } \right)^{\frac{{Q - p}}{{pQ}}}} \nonumber \\
   &  \le c\left( {n,p} \right){\left| {{B_{6r}}} \right|^{1 - \frac{{Q - p}}{{pQ}}}}{\left( {\int_{{B_{6r}}} {{{\left| {{\nabla _H}\upsilon } \right|}^p}d\xi } } \right)^{\frac{1}{p}}}\nonumber \\
  & \le c\left( {n,p} \right){\left| {{B_{6r}}} \right|^{1 - \frac{{Q - p}}{{pQ}}}}{\left( {\int_{{B_{6r}}} {{{\left| {{\nabla _H}\log \omega } \right|}^p}d\xi } } \right)^{\frac{1}{p}}} \nonumber \\
   & \le c\left( {n,p} \right){\left| {{B_{6r}}} \right|^{1 - \frac{{Q - p}}{{pQ}}}}{r^{\frac{{Q - p}}{p}}} \le c\left| {{B_{6r}}} \right|.
\end{align}
Note the relation
\[\left\{ {\upsilon  = 0} \right\} = \left\{ {\omega  \ge k + {d_\varepsilon }} \right\} = \left\{ {u \ge k} \right\},\]
we see by the assumption \eqref{eq56} that
\begin{equation}\label{eq515}
\left| {{B_{6r}}\left( {{\xi _0}} \right) \cap \left\{ {\upsilon  = 0} \right\}} \right| \ge \frac{\tau }{{{6^Q}}}\left| {{B_{6r}}\left( {{\xi _0}} \right)} \right|
\end{equation}
and thus
\begin{align}\label{eq516}
   \log \frac{1}{{2\delta }}& = \fint_{{B_{6r}}\left( {{\xi _0}} \right) \cap \left\{ {\upsilon  = 0} \right\}} {\left( {\log \frac{1}{{2\delta }} - \upsilon \left( \xi  \right)} \right)d\xi }  \nonumber\\
   & \le \frac{{{6^Q}}}{\tau }\left( {\log \frac{1}{{2\delta }} - {{\left( \upsilon  \right)}_{{B_{6r}}}}} \right).
\end{align}
Integrating both sides of the above inequality over ${B_{6r}}\left( {{\xi _0}} \right) \cap \left\{ {\upsilon  = \log \frac{1}{{2\delta }}} \right\}$ and using \eqref{eq514}, we gain
\begin{align*}
   \left| {{B_{6r}}\left( {{\xi _0}} \right) \cap \left\{ {\upsilon  = \log \frac{1}{{2\delta }}} \right\}} \right|\log \frac{1}{{2\delta }}& \le \frac{{{6^Q}}}{\tau }\int_{{B_{6r}}\left( {{\xi _0}} \right) \cap \left\{ {\upsilon  = \log \frac{1}{{2\delta }}} \right\}} {\left( {\log \frac{1}{{2\delta }} - {{\left( \upsilon  \right)}_{{B_{6r}}}}} \right)d\xi }  \\
   &  \le \frac{{{6^Q}}}{\tau }\int_{{B_{6r}}\left( {{\xi _0}} \right)} {\left( {\upsilon \left( \xi  \right) - {{\left( \upsilon  \right)}_{{B_{6r}}}}} \right)d\xi }  \le \frac{c}{\tau }\left| {{B_{6r}}\left( {{\xi _0}} \right)} \right|.
\end{align*}
Therefore, for any $\delta  \in \left( {0,\frac{1}{4}} \right)$, we have
\[\left| {{B_{6r}}\left( {{\xi _0}} \right) \cap \left\{ {\omega  \le 2\delta \left( {k + {d_\varepsilon }} \right)} \right\}} \right| \le \frac{c}{\tau }\frac{1}{{\log \frac{1}{{2\delta }}}}\left| {{B_{6r}}\left( {{\xi _0}} \right)} \right|,\]
which means that \eqref{eq58} holds.

\textbf{Step 2.} For any $\varepsilon  > 0$, we prove that there exists a constant $\delta  = \delta \left( {n,p,s,\tau } \right) \in \left( {0,\frac{1}{4}} \right)$ such that
\begin{equation}\label{eq517}
\mathop {{\rm{ess inf}}}\limits_{{B_{4r}}\left( {{\xi _0}} \right)} u \ge \delta k - {\left( {\frac{r}{R}} \right)^{\frac{p}{{p - 1}}}}{\rm{Tail}}({u_ - };{\xi _0},R) - 2\varepsilon .
\end{equation}
Then it follows \eqref{eq57} by the arbitrariness of $\varepsilon $ and \eqref{eq517}.

To prove \eqref{eq517}, without loss of generality, we assume that
\begin{equation}\label{eq518}
\delta k \ge {\left( {\frac{r}{R}} \right)^{\frac{p}{{p - 1}}}}{\rm{Tail}}({u_ - };{\xi _0},R) + 2\varepsilon .
\end{equation}
Otherwise, \eqref{eq517} holds trivially because $u \ge 0$ in ${B_R}\left( {{\xi _0}} \right)$.

Let $\rho  \in \left[ {r,6r} \right]$ and $\psi  \in C_0^\infty \left( {{B_\rho }\left( {{\xi _0}} \right)} \right)$ be a cut-off function satisfying $0 \le \psi  \le 1$ in ${B_\rho }\left( {{\xi _0}} \right)$. Denote $\omega  = {\left( {l - u} \right)_ + },\;l \in \left( {\delta k,2\delta k} \right)$, and take $\phi  = \omega {\psi ^p}$ as the test function in \eqref{eq22}. Similar to the proof of Lemma \ref{Le27}, we obtain
\begin{align}\label{eq519}
   & \int_{{B_\rho }} {{\psi ^p}{{\left| {{\nabla _H}\omega } \right|}^p}d\xi }  + \int_{{B_\rho }} {\int_{{B_\rho }} {\frac{{{{\left| {\omega \left( \xi  \right)\psi \left( \xi  \right) - \omega \left( \eta  \right)\psi \left( \eta  \right)} \right|}^p}}}{{\left\| {{\eta ^{ - 1}} \circ \xi } \right\|_{{{\rm{H}}^n}}^{Q + sp}}}{\mkern 1mu} d\xi d\eta } }  \nonumber\\
   \le& c(\int_{{B_\rho }} {{\omega ^p}{{\left| {{\nabla _H}\psi } \right|}^p}d\xi }  + \int_{{B_\rho }} {\int_{{B_\rho }} {\max {{\{ \omega \left( \xi  \right),\omega \left( \eta  \right)\} }^p}\frac{{{{\left| {\psi \left( \xi  \right) - \psi \left( \eta  \right)} \right|}^p}}}{{\left\| {{\eta ^{ - 1}} \circ \xi } \right\|_{{{\rm{H}}^n}}^{Q + sp}}}{\mkern 1mu} d\xi d\eta } }\nonumber\\
   & + \mathop {{\rm{ess}}\;\sup }\limits_{\xi  \in {\rm{supp}}\;\psi } \int_{{{\rm{H}}^n}\backslash {B_\rho }} {\frac{{\omega {{\left( \eta  \right)}^{p - 1}}}}{{\left\| {{\eta ^{ - 1}} \circ \xi } \right\|_{{{\rm{H}}^n}}^{Q + sp}}}d\eta }  \cdot \int_{{B_\rho }} {\omega {\psi ^p}d\xi } ).
\end{align}
Furthermore, by using the estimates of ${J_2} + {J_3}$ in Lemma 4.2 of \cite{PP22}, we get
\begin{align}\label{eq520}
  & \int_{{B_\rho }} {{\psi ^p}{{\left| {{\nabla _H}\omega } \right|}^p}d\xi }  + \int_{{B_\rho }} {\int_{{B_\rho }} {\frac{{{{\left| {\omega \left( \xi  \right)\psi \left( \xi  \right) - \omega \left( \eta  \right)\psi \left( \eta  \right)} \right|}^p}}}{{\left\| {{\eta ^{ - 1}} \circ \xi } \right\|_{{{\rm{H}}^n}}^{Q + sp}}}{\mkern 1mu} d\xi d\eta } }  \nonumber\\
  \le &  c\int_{{B_\rho }} {{\omega ^p}{{\left| {{\nabla _H}\psi } \right|}^p}d\xi }  + c\int_{{B_\rho }} {\int_{{B_\rho }} {\max {{\{ \omega \left( \xi  \right),\omega \left( \eta  \right)\} }^p}\frac{{{{\left| {\psi \left( \xi  \right) - \psi \left( \eta  \right)} \right|}^p}}}{{\left\| {{\eta ^{ - 1}} \circ \xi } \right\|_{{{\rm{H}}^n}}^{Q + sp}}}{\mkern 1mu} d\xi d\eta } }\nonumber\\
 & + cl\mathop {{\rm{ess}}\;\sup }\limits_{\xi  \in {\rm{supp}}\;\psi } \int_{{{\rm{H}}^n}\backslash {B_\rho }} {\frac{{{{\left( {l + {{\left( {u\left( \eta  \right)} \right)}_ - }} \right)}^{p - 1}}}}{{\left\| {{\eta ^{ - 1}} \circ \xi } \right\|_{{{\rm{H}}^n}}^{Q + sp}}}d\eta }  \cdot \left| {{B_\rho }\left( {{\xi _0}} \right) \cap \left\{ {u < l} \right\}} \right| \nonumber\\
  = :&{I_1} + {I_2} + {I_3}.
\end{align}

For $j = 0,1,2, \cdots $, denote
\begin{equation}\label{eq521}
l = {k_j} = \delta k + {2^{ - j - 1}}\delta k,\;\;\rho  = {\rho _j} = 4r + {2^{1 - j}}r,\;\;{\bar \rho _j} = \frac{{{\rho _j} + {\rho _{j + 1}}}}{2}.
\end{equation}
Then, for any $j = 0,1,2, \cdots $, we have $l \in \left( {\delta k,2\delta k} \right),\;{\rho _j},{\bar \rho _j} \in \left( {4r,6r} \right)$ and
\begin{equation}\label{eq522}
{k_j} - {k_{j + 1}} = {2^{ - j - 2}}\delta k \ge {2^{ - j - 3}}{k_j}.
\end{equation}
Let ${B_j} = {B_{{\rho _j}}}\left( {{\xi _0}} \right),\;{\bar B_j} = {B_{{{\bar \rho }_j}}}\left( {{\xi _0}} \right)$. Note
\[{\omega _j}: = {\left( {{k_j} - u} \right)_ + } \ge {2^{ - j - 3}}{k_j}{\chi _{\left\{ {u < {k_{j + 1}}} \right\}}}.\]
Let $\left( {{\psi _j}} \right)_{j = 0}^\infty  \subset C_0^\infty \left( {{{\bar B}_j}} \right),j = 0,1,2, \cdots $ be a sequence of cut-off functions with
\begin{center}
$0 \le {\psi _j} \le 1,\;\left| {{\nabla _H}{\psi _j}} \right| \le \frac{{{2^{j + 3}}}}{r}$ in ${\bar B_j}$ and ${\psi _j} = 1$ in ${B_{j + 1}}$.
\end{center}
We choose $\psi  = {\psi _j}$ and $\omega  = {\omega _j}$ in \eqref{eq520} and obtain by the properties of ${\psi _j}$ that
\begin{equation}\label{eq523}
{I_1} = c\int_{{B_j}} {\omega _j^p{{\left| {{\nabla _H}{\psi _j}} \right|}^p}d\xi }  \le c\left( p \right){2^{jp}}k_j^p{r^{ - p}}\left| {{B_j} \cap \left\{ {u < {k_j}} \right\}} \right|.
\end{equation}
By Lemma 4.2 in [\cite{PP22}, p17], we have
\begin{align}\label{eq524}
  {I_2} & = c\int_{{B_j}} {\int_{{B_j}} {\max {{\{ {\omega _j}\left( \xi  \right),{\omega _j}\left( \eta  \right)\} }^p}\frac{{{{\left| {{\psi _j}\left( \xi  \right) - {\psi _j}\left( \eta  \right)} \right|}^p}}}{{\left\| {{\eta ^{ - 1}} \circ \xi } \right\|_{{{\rm{H}}^n}}^{Q + sp}}}{\mkern 1mu} d\xi d\eta } }  \nonumber\\
   &  \le c\left( {n,p,s} \right){2^{jp}}k_j^p{r^{ - sp}}\left| {{B_j} \cap \left\{ {u < {k_j}} \right\}} \right| \nonumber\\
   & \le c\left( {n,p,s} \right){2^{jp}}k_j^p{r^{ - p}}\left| {{B_j} \cap \left\{ {u < {k_j}} \right\}} \right|.
\end{align}
By (4.18) in [\cite{PP22}, Lemma 4.2], it follows
\begin{align}\label{eq525}
   {I_3}& = c{k_j}\mathop {{\rm{ess}}\;\sup }\limits_{\xi  \in {\rm{supp}}\;{\psi _j}} \int_{{{\rm{H}}^n}\backslash {B_j}} {\frac{{{{\left( {{k_j} + {{\left( {u\left( \eta  \right)} \right)}_ - }} \right)}^{p - 1}}}}{{\left\| {{\eta ^{ - 1}} \circ \xi } \right\|_{{{\rm{H}}^n}}^{Q + sp}}}d\eta }  \cdot \left| {{B_j}\left( {{\xi _0}} \right) \cap \left\{ {u < {k_j}} \right\}} \right| \nonumber\\
   &  \le c\left( {n,p,s} \right){2^{j\left( {Q + sp} \right)}}k_j^p{r^{ - sp}}\left| {{B_j} \cap \left\{ {u < {k_j}} \right\}} \right|\nonumber\\
   & \le c\left( {n,p,s} \right){2^{j\left( {Q + sp} \right)}}k_j^p{r^{ - p}}\left| {{B_j} \cap \left\{ {u < {k_j}} \right\}} \right|.
\end{align}
Using \eqref{eq520} and \eqref{eq523}-\eqref{eq525}, we see
\begin{equation}\label{eq526}
\int_{{B_j}} {\psi _j^p{{\left| {{\nabla _H}{\omega _j}} \right|}^p}d\xi }  \le c\left( {n,p,s} \right){2^{j\left( {Q + sp + p} \right)}}k_j^p{r^{ - p}}\left| {{B_j} \cap \left\{ {u < {k_j}} \right\}} \right|.
\end{equation}
Therefore, by Lemma \ref{Le22}, \eqref{eq523} and \eqref{eq526}, we deduce
\begin{align}\label{eq527}
   & {\left( {{k_j} - {k_{j + 1}}} \right)^p}{\left( {\frac{{\left| {{B_{j + 1}} \cap \left\{ {u < {k_{j + 1}}} \right\}} \right|}}{{\left| {{B_{j + 1}}} \right|}}} \right)^{^{\frac{p}{{p*}}}}} \le {\left( {\fint_{{B_{j + 1}}} {\omega _j^{p*}\psi _j^{p*}d\xi } } \right)^{^{\frac{p}{{p*}}}}}\nonumber \\
  \le &  c{\left( {\fint_{{B_j}} {\omega _j^{p*}\psi _j^{p*}d\xi } } \right)^{^{\frac{p}{{p*}}}}} \le c{r^p}\fint_{{B_j}} {{{\left| {{\nabla _H}\left( {{\omega _j}{\psi _j}} \right)} \right|}^p}d\xi }\nonumber \\
   \le& c\left( {n,p,s} \right){2^{j\left( {Q + sp + p} \right)}}k_j^p\frac{{\left| {{B_j} \cap \left\{ {u < {k_j}} \right\}} \right|}}{{\left| {{B_j}} \right|}}.
\end{align}

Let
\[{Y_j} = \frac{{\left| {{B_j} \cap \left\{ {u < {k_j}} \right\}} \right|}}{{\left| {{B_j}} \right|}},\;j = 0,1,2, \cdots .\]
It gives from \eqref{eq527} and \eqref{eq522} that
\begin{equation}\label{eq528}
{Y_{j + 1}} \le c\left( {n,p,s} \right){2^{j\left( {Q + sp + 2p} \right)\frac{{p*}}{p}}}Y_j^{\frac{{p*}}{p}},\;j = 0,1,2, \cdots .
\end{equation}
Taking ${c_0} = c\left( {n,p,s} \right),\;b = {2^{\left( {Q + sp + 2p} \right)\frac{{p*}}{p}}},\;\beta  = \frac{{p*}}{p} - 1$ as in Lemma \ref{Le26}, it follows from \eqref{eq518} that
\[{k_0} = \frac{3}{2}\delta k \le 2\delta k - \frac{1}{2}{\left( {\frac{r}{R}} \right)^{\frac{p}{{p - 1}}}}{\rm{Tail}}({u_ - };{\xi _0},R) - \varepsilon .\]
By \eqref{eq58}, for each $\delta  \in \left( {0,\frac{1}{4}} \right)$ and $\varepsilon  > 0$, we have
\begin{equation}\label{eq529}
{Y_0} \le \frac{{\left| {{B_{6r}}\left( {{\xi _0}} \right) \cap \left\{ {u \le 2\delta k - \frac{1}{2}{{\left( {\frac{r}{R}} \right)}^{\frac{{sp}}{{p - 1}}}}{\rm{Tail}}({u_ - };{\xi _0},R) - \varepsilon } \right\}} \right|}}{{\left| {{B_{6r}}\left( {{\xi _0}} \right)} \right|}} \le \frac{{{c_1}}}{{\tau \log \frac{1}{{2\delta }}}}.
\end{equation}
Choose $\delta  = \delta \left( {n,p,s,\tau } \right) \in \left( {0,\frac{1}{4}} \right)$ such that
\[0 < \delta  = \frac{1}{4}{e^{ - \frac{{{c_1}c_0^{\frac{1}{\beta }}{b^{\frac{1}{{{\beta ^2}}}}}}}{\tau }}} < \frac{1}{4}.\]
Then we have by combining with \eqref{eq529} that ${Y_0} \le c_0^{ - \frac{1}{\beta }}{b^{ - \frac{1}{{{\beta ^2}}}}}.$ Thus, applying Lemma \ref{Le26} implies $\mathop {\lim }\limits_{j \to \infty } {Y_j} = 0$, i.e.
\[\mathop {{\rm{ess inf}}}\limits_{{B_{4r}}\left( {{\xi _0}} \right)} u \ge \delta k.\]
Hence, \eqref{eq517} holds.

Using the arbitrariness of $\varepsilon $ in \eqref{eq517}, we prove \eqref{eq57}.
\end{proof}

Similar to the proof of Lemma 4.1 in \cite{DKP14}, we have the following lemma by using the classical Krylov-Safonov covering lemma and Lemma \ref{Le52}.

\begin{lemma}\label{Le53}
Let $u \in HW_{loc}^{1,p}\left( \Omega  \right)\;\left( {1 < p < Q} \right)$ be a weak supersolution to \eqref{eq15} and $u \ge 0$ in ${B_R}\left( {{\xi _0}} \right) \subset \Omega $. Then there exist constants $\lambda  = \lambda \left( {n,p,s} \right) \in \left( {0,1} \right)$ and $c = c\left( {n,p,s} \right) \ge 1$ such that
\begin{equation}\label{eq530}
{\left( {\fint_{{B_r}\left( {{\xi _0}} \right)} {{u^\lambda }d\xi } } \right)^{\frac{1}{\lambda }}} \le c\mathop {{\rm{ess inf}}}\limits_{{B_r}\left( {{\xi _0}} \right)} u + c{\left( {\frac{r}{R}} \right)^{\frac{p}{{p - 1}}}}{\rm{Tail}}({u_ - };{\xi _0},R),
\end{equation}
where $r \in \left( {0,1} \right]$ and ${B_r}\left( {{\xi _0}} \right) \subset {B_R}\left( {{\xi _0}} \right)$. Here ${\rm{Tail}}\left(  \cdot  \right)$ is given by \eqref{eq21}.
\end{lemma}

\begin{lemma}[Iteration lemma, \cite{GG82}]\label{Le54}
For $0 \le {T_0} \le t \le {T_1},$ let $f:\left[ {{T_0},{T_1}} \right] \to \left[ {0,\infty } \right)$ be a non-negative bounded function. If for all ${T_0} \le t < s \le {T_1}$, it holds
\begin{equation}\label{eq531}
f\left( t \right) \le A{\left( {s - t} \right)^{ - \alpha }} + B + \theta f\left( s \right),
\end{equation}
where $A,B,\alpha ,\theta  < 1$ are non-negative constants, then there exists a constant $c = c\left( {\alpha ,\theta } \right)$ such that for any ${T_0} \le \rho  < R \le {T_1}$,
\begin{equation}\label{eq532}
f\left( \rho  \right) \le c\left( {A{{\left( {R - \rho } \right)}^{ - \alpha }} + B} \right).
\end{equation}
\end{lemma}

\textbf{Proof of Theorem \ref{Th13}.} Let $0 < \rho  < r$. Since $u$ is a weak solution to \eqref{eq15}, we see that $u$ is both a weak supersolution and a weak subsolution to \eqref{eq15}. Thus, by Theorem \ref{Th11}, we have
\begin{equation}\label{eq533}
\mathop {{\rm{ess}}\sup u}\limits_{{B_{\frac{\rho }{2}}}\left( {{\xi _0}} \right)}  \le \delta {\rm{Tail}}({u_ + };{\xi _0},\frac{\rho }{2}) + c{\delta ^{ - \frac{{\left( {p - 1} \right)Q}}{{{p^2}}}}}{\left( {\fint_{{B_\rho }} {u_ + ^pd\xi } } \right)^{\frac{1}{p}}}.
\end{equation}
Furthermore, we obtain by using Lemma \ref{Le51} and \eqref{eq533} that
\begin{equation}\label{eq534}
\mathop {{\rm{ess}}\sup u}\limits_{{B_{\frac{\rho }{2}}}\left( {{\xi _0}} \right)}  \le c\delta \left( {\mathop {{\rm{ess}}\sup }\limits_{{B_\rho }\left( {{\xi _0}} \right)} u + {{\left( {\frac{\rho }{R}} \right)}^{\frac{p}{{p - 1}}}}{\rm{Tail}}({u_ - };{\xi _0},R)} \right) + c{\delta ^{ - \frac{{\left( {p - 1} \right)Q}}{{{p^2}}}}}{\left( {\fint_{{B_\rho }} {u_ + ^pd\xi } } \right)^{\frac{1}{p}}},
\end{equation}
where $c = c\left( {n,p,s} \right)$. Let $\frac{1}{2} \le \sigma ' < \sigma  \le 1,\;\rho  = \left( {\sigma  - \sigma '} \right)r$. For any $\lambda  \in \left( {0,p} \right)$, we have by means of the covering argument that
\begin{align}\label{eq535}
\mathop {{\rm{ess}}\sup }\limits_{{B_{\sigma 'r}}\left( {{\xi _0}} \right)} u \le & c\frac{{{\delta ^{ - \frac{{\left( {p - 1} \right)Q}}{{{p^2}}}}}}}{{{{\left( {\sigma  - \sigma '} \right)}^{\frac{Q}{p}}}}}{\left( {\fint_{{B_{\sigma r}}} {{u^p}d\xi } } \right)^{\frac{1}{p}}} + c\delta \mathop {{\rm{ess}}\sup }\limits_{{B_{\sigma r}}\left( {{\xi _0}} \right)} u + c\delta {\left( {\frac{r}{R}} \right)^{\frac{p}{{p - 1}}}}{\rm{Tail}}({u_ - };{\xi _0},R)\nonumber \\
   =& c\frac{{{\delta ^{ - \frac{{\left( {p - 1} \right)Q}}{{{p^2}}}}}}}{{{{\left( {\sigma  - \sigma '} \right)}^{\frac{Q}{p}}}}}{\left( {\fint_{{B_{\sigma r}}} {{u^{p - \lambda }} \cdot {u^\lambda }d\xi } } \right)^{\frac{1}{p}}} + c\delta \mathop {{\rm{ess}}\sup }\limits_{{B_{\sigma r}}\left( {{\xi _0}} \right)} u + c\delta {\left( {\frac{r}{R}} \right)^{\frac{p}{{p - 1}}}}{\rm{Tail}}({u_ - };{\xi _0},R)\nonumber \\
\le& c\frac{{{\delta ^{ - \frac{{\left( {p - 1} \right)Q}}{{{p^2}}}}}}}{{{{\left( {\sigma  - \sigma '} \right)}^{\frac{Q}{p}}}}}{\left( {\mathop {{\rm{ess}}\sup }\limits_{{B_{\sigma r}}\left( {{\xi _0}} \right)} u} \right)^{\frac{{p - \lambda }}{p}}}{\left( {\fint_{{B_{\sigma r}}} {{u^\lambda }d\xi } } \right)^{\frac{1}{p}}} \nonumber \\
    &+ c\delta \mathop {{\rm{ess}}\sup }\limits_{{B_{\sigma r}}\left( {{\xi _0}} \right)} u + c\delta {\left( {\frac{r}{R}} \right)^{\frac{p}{{p - 1}}}}{\rm{Tail}}({u_ - };{\xi _0},R).
\end{align}
For any $\lambda  \in \left( {0,p} \right)$, using Young's inequality with exponents $\frac{p}{\lambda }$ and $\frac{p}{{p - \lambda }}$, and taking $\delta  = \frac{1}{{4c}}$ in \eqref{eq535}, we get
\begin{equation}\label{eq536}
\mathop {{\rm{ess}}\sup }\limits_{{B_{\sigma 'r}}\left( {{\xi _0}} \right)} u \le \frac{1}{2}\mathop {{\rm{ess}}\sup }\limits_{{B_{\sigma r}}\left( {{\xi _0}} \right)} u + \frac{c}{{{{\left( {\sigma  - \sigma '} \right)}^{\frac{Q}{\lambda }}}}}{\left( {\fint_{{B_{\sigma r}}} {{u^\lambda }d\xi } } \right)^{\frac{1}{\lambda }}} + c{\left( {\frac{r}{R}} \right)^{\frac{p}{{p - 1}}}}{\rm{Tail}}({u_ - };{\xi _0},R).
\end{equation}
Thus by Lemma \ref{Le54}, it follows that for any $\lambda  \in \left( {0,p} \right)$,
\begin{equation}\label{eq537}
\mathop {{\rm{ess}}\sup u}\limits_{{B_{\frac{r}{2}}}\left( {{\xi _0}} \right)}  \le c{\left( {\fint_{{B_r}} {{u^\lambda }d\xi } } \right)^{\frac{1}{\lambda }}} + c{\left( {\frac{r}{R}} \right)^{\frac{p}{{p - 1}}}}{\rm{Tail}}({u_ - };{\xi _0},R).
\end{equation}
Therefore, we conclude \eqref{eq18} by combining Lemma \ref{Le53} and \eqref{eq537}.

\section{Weak Harnack Inequality}
\label{Section 6}

In this section, we first state a known Lemma \ref{Le61}, then use it to derive a key lemma, namely Lemma \ref{Le62}, and finally use key lemma together with Lemma \ref{Le53} to prove Theorem \ref{Th14}.

\begin{lemma}[\cite{DKP16}]\label{Le61}
Let $p \ge 1,\;\varepsilon  \in \left( {0,1} \right]$, then for any $a,b \in {\mathbb{R}^m}\;(m \ge 1)$, we have
\[{\left| a \right|^p} \le \left( {1 + {c_p}\varepsilon } \right){\left| b \right|^p} + \left( {1 + {c_p}\varepsilon } \right){\varepsilon ^{1 - p}}{\left| {a - b} \right|^p},\]
where ${c_p}: = \left( {p - 1} \right)\Gamma \left( {\max \{ 1,p - 2\} } \right)$ and $\Gamma$ denotes the standard Gamma function.
\end{lemma}

\begin{lemma}\label{Le62}
For $1 < q < p,\;d > 0,$ let $u \in HW_{loc}^{1,p}\left( \Omega  \right)$ be a weak supersolution to \eqref{eq15} satisfying $u \ge 0$ in ${B_R}\left( {{\xi _0}} \right) \subset \Omega $. Denote $\omega  = {\left( {u + d} \right)^{\frac{{p - q}}{p}}}$, then for ${B_r} \equiv {B_r}\left( {{\xi _0}} \right) \subset {B_{\frac{{3R}}{4}}}\left( {{\xi _0}} \right)$ and non-negative functions $\psi  \in C_0^\infty \left( {{B_r}\left( {{\xi _0}} \right)} \right)$, there exists a positive constant $c = c\left( {n,p,s} \right)$ such that
\begin{align}\label{eq61}
   \int_{{B_r}} {{{\left| {{\nabla _H}\omega } \right|}^p}{\psi ^p}d\xi }  \le& c(\frac{{{{\left( {p - q} \right)}^p}}}{{{{\left( {q - 1} \right)}^{\frac{p}{{p - 1}}}}}}\int_{{B_r}} {{\omega ^p}{{\left| {{\nabla _H}\psi } \right|}^p}d\xi } \nonumber \\
   &  + \frac{{{{\left( {p - q} \right)}^p}}}{{{{\left( {q - 1} \right)}^p}}}\int_{{B_r}} {\int_{{B_r}} {\frac{{{{\left( {\max \left\{ {\omega \left( \xi  \right),\omega \left( \eta  \right)} \right\}} \right)}^p}{{\left| {\psi \left( \xi  \right) - \psi \left( \eta  \right)} \right|}^p}}}{{\left\| {{\eta ^{ - 1}} \circ \xi } \right\|_{{{\rm{H}}^n}}^{Q + sp}}}d\xi d\eta } }\nonumber \\
& + \frac{{{{\left( {p - q} \right)}^p}}}{{q - 1}}(\mathop {{\rm{ess}}\sup }\limits_{\xi  \in {\rm{supp}}\;\psi } \int_{{{\rm{H}}^n}\backslash {B_r}} {\frac{1}{{\left\| {{\eta ^{ - 1}} \circ \xi } \right\|_{{{\rm{H}}^n}}^{Q + sp}}}d\eta } \nonumber \\
& + {d^{1 - p}}{R^{ - p}}{\rm{Tail}}{({u_ - };{\xi _0},R)^{p - 1}}\int_{{B_r}} {\omega {{\left( \xi  \right)}^p}\psi {{\left( \xi  \right)}^p}d\xi } ),
\end{align}
where ${\rm{Tail}}\left(  \cdot  \right)$ is given by \eqref{eq21}.
\end{lemma}

\begin{proof}
Let $d > 0,\;v = u + d$, then $v \in HW_{loc}^{1,p}\left( \Omega  \right)$ is a weak supersolution to \eqref{eq15}. For $1 < q < p$, we can take some sufficiently small $\varepsilon  > 0$, such that $q \in \left[ {1 + \varepsilon ,p - \varepsilon } \right]$. Choosing $\phi  = {v^{1 - q}}{\psi ^p}$ as the test function in \eqref{eq22}, we obtain
\begin{align}\label{eq62}
 0 \le  & \int_{{B_r}} {{{\left| {{\nabla _H}u} \right|}^{p - 2}}{\nabla _H}u \cdot {\nabla _H}\left( {{v^{1 - q}}{\psi ^p}} \right)d\xi } \nonumber \\
   &  + \int_{{{\rm{H}}^n}} {\int_{{{\rm{H}}^n}} {\frac{{{{\left| {u\left( \xi  \right) - u\left( \eta  \right)} \right|}^{p - 2}}\left( {u\left( \xi  \right) - u\left( \eta  \right)} \right)\left( {v{{\left( \xi  \right)}^{1 - q}}\psi {{\left( \xi  \right)}^p} - v{{\left( \eta  \right)}^{1 - q}}\psi {{\left( \eta  \right)}^p}} \right)}}{{\left\| {{\eta ^{ - 1}} \circ \xi } \right\|_{{{\rm{H}}^n}}^{Q + sp}}}d\xi d\eta } }\nonumber \\
   = &\int_{{B_r}} {{{\left| {{\nabla _H}u} \right|}^{p - 2}}{\nabla _H}u \cdot {\nabla _H}\left( {{v^{1 - q}}{\psi ^p}} \right)d\xi }  \nonumber \\
   & + \int_{{B_r}} {\int_{{B_r}} {\frac{{{{\left| {u\left( \xi  \right) - u\left( \eta  \right)} \right|}^{p - 2}}\left( {u\left( \xi  \right) - u\left( \eta  \right)} \right)\left( {v{{\left( \xi  \right)}^{1 - q}}\psi {{\left( \xi  \right)}^p} - v{{\left( \eta  \right)}^{1 - q}}\psi {{\left( \eta  \right)}^p}} \right)}}{{\left\| {{\eta ^{ - 1}} \circ \xi } \right\|_{{{\rm{H}}^n}}^{Q + sp}}}d\xi d\eta } } \nonumber \\
   & + 2\int_{{{\rm{H}}^n}\backslash {B_r}} {\int_{{B_r}} {\frac{{{{\left| {u\left( \xi  \right) - u\left( \eta  \right)} \right|}^{p - 2}}\left( {u\left( \xi  \right) - u\left( \eta  \right)} \right)v{{\left( \xi  \right)}^{1 - q}}\psi {{\left( \xi  \right)}^p}}}{{\left\| {{\eta ^{ - 1}} \circ \xi } \right\|_{{{\rm{H}}^n}}^{Q + sp}}}d\xi d\eta } } \nonumber \\
    =& \int_{{B_r}} {{{\left| {{\nabla _H}v} \right|}^{p - 2}}{\nabla _H}v \cdot {\nabla _H}\left( {{v^{1 - q}}{\psi ^p}} \right)d\xi } \nonumber \\
  & + \int_{{B_r}} {\int_{{B_r}} {\frac{{{{\left| {v\left( \xi  \right) - v\left( \eta  \right)} \right|}^{p - 2}}\left( {v\left( \xi  \right) - v\left( \eta  \right)} \right)\left( {v{{\left( \xi  \right)}^{1 - q}}\psi {{\left( \xi  \right)}^p} - v{{\left( \eta  \right)}^{1 - q}}\psi {{\left( \eta  \right)}^p}} \right)}}{{\left\| {{\eta ^{ - 1}} \circ \xi } \right\|_{{{\rm{H}}^n}}^{Q + sp}}}d\xi d\eta } }  \nonumber \\
  & + 2\int_{{{\rm{H}}^n}\backslash {B_r}} {\int_{{B_r}} {\frac{{{{\left| {v\left( \xi  \right) - v\left( \eta  \right)} \right|}^{p - 2}}\left( {v\left( \xi  \right) - v\left( \eta  \right)} \right)v{{\left( \xi  \right)}^{1 - q}}\psi {{\left( \xi  \right)}^p}}}{{\left\| {{\eta ^{ - 1}} \circ \xi } \right\|_{{{\rm{H}}^n}}^{Q + sp}}}d\xi d\eta } }  \nonumber \\
   = :&{I_1} + {I_2} + 2{I_3}.
\end{align}

For ${I_1}$, it follows from Young's inequality that
\begin{align}\label{eq63}
   {I_1} = & \int_{{B_r}} {{{\left| {{\nabla _H}v} \right|}^{p - 2}}{\nabla _H}v \cdot {\nabla _H}\left( {{v^{1 - q}}{\psi ^p}} \right)d\xi }   \nonumber \\
  = &  \left( {1 - q} \right)\int_{{B_r}} {{{\left| {{\nabla _H}v} \right|}^{p - 2}}{\nabla _H}v \cdot {\psi ^p}{v^{ - q}}{\nabla _H}vd\xi } \nonumber \\
   & + p\int_{{B_r}} {{{\left| {{\nabla _H}v} \right|}^{p - 2}}{\nabla _H}v \cdot {v^{1 - q}}{\psi ^{p - 1}}{\nabla _H}\psi d\xi } \nonumber \\
    \le& \left( {1 - q} \right)\int_{{B_r}} {{v^{ - q}}{{\left| {{\nabla _H}v} \right|}^p}{\psi ^p}d\xi }  \nonumber \\
     & + p\int_{{B_r}} {\left( {{v^{\frac{q}{p} - q}}{{\left| {{\nabla _H}v} \right|}^{p - 1}}{\psi ^{p - 1}}} \right) \cdot \left( {{v^{1 - \frac{q}{p}}}\left| {{\nabla _H}\psi } \right|} \right)d\xi } \nonumber \\
       \le& \left( {1 - q} \right)\int_{{B_r}} {{v^{ - q}}{{\left| {{\nabla _H}v} \right|}^p}{\psi ^p}d\xi } \nonumber \\
       & + \frac{{q - 1}}{2}\int_{{B_r}} {{v^{ - q}}{{\left| {{\nabla _H}v} \right|}^p}{\psi ^p}d\xi }  + \frac{{c\left( p \right)}}{{{{\left( {q - 1} \right)}^{\frac{1}{{p - 1}}}}}}\int_{{B_r}} {{v^{p - q}}{{\left| {{\nabla _H}\psi } \right|}^p}d\xi } \nonumber \\
         = &\frac{{1 - q}}{2}\int_{{B_r}} {{v^{ - q}}{{\left| {{\nabla _H}v} \right|}^p}{\psi ^p}d\xi }  + \frac{{c\left( p \right)}}{{{{\left( {q - 1} \right)}^{\frac{1}{{p - 1}}}}}}\int_{{B_r}} {{v^{p - q}}{{\left| {{\nabla _H}\psi } \right|}^p}d\xi } \nonumber \\
          = & - \frac{{q - 1}}{2}{\left( {\frac{p}{{p - q}}} \right)^p}\int_{{B_r}} {{{\left| {{\nabla _H}\left( {{v^{\frac{{p - q}}{p}}}} \right)} \right|}^p}{\psi ^p}d\xi }  + \frac{{c\left( p \right)}}{{{{\left( {q - 1} \right)}^{\frac{1}{{p - 1}}}}}}\int_{{B_r}} {{v^{p - q}}{{\left| {{\nabla _H}\psi } \right|}^p}d\xi } .
\end{align}

For ${I_3}$, we note that for any $\xi  \in {B_r} \subset {B_R}$ and any $\eta  \in {{\rm{H}}^n}{\rm{\backslash }}{B_r}$, it holds
\begin{align*}
   {\left| {v\left( \xi  \right) - v\left( \eta  \right)} \right|^{p - 2}}\left( {v\left( \xi  \right) - v\left( \eta  \right)} \right)& = {\left| {u\left( \xi  \right) - u\left( \eta  \right)} \right|^{p - 2}}\left( {u\left( \xi  \right) - u\left( \eta  \right)} \right) \\
   &  \le cu{\left( \xi  \right)^{p - 1}} - c\left( {u\left( \eta  \right)} \right)_ + ^{p - 1} + c\left( {u\left( \eta  \right)} \right)_ - ^{p - 1}\\
   & \le cv{\left( \xi  \right)^{p - 1}} + c\left( {u\left( \eta  \right)} \right)_ - ^{p - 1}
\end{align*}
and
\[v{\left( \xi  \right)^{1 - q}} \le {d^{1 - p}}v{\left( \xi  \right)^{p - q}}.\]
Moreover, from the definition of $u$, it follows ${\left( {u\left( \eta  \right)} \right)_ - } = 0$ for any $\eta  \in {B_R}$.
Therefore,
\begin{align*}
   {I_3}& = \int_{{{\rm{H}}^n}\backslash {B_r}} {\int_{{B_r}} {\frac{{{{\left| {v\left( \xi  \right) - v\left( \eta  \right)} \right|}^{p - 2}}\left( {v\left( \xi  \right) - v\left( \eta  \right)} \right)v{{\left( \xi  \right)}^{1 - q}}\psi {{\left( \xi  \right)}^p}}}{{\left\| {{\eta ^{ - 1}} \circ \xi } \right\|_{{{\rm{H}}^n}}^{Q + sp}}}d\xi d\eta } }  \\
   &  \le c\left( p \right)\left( {\int_{{{\rm{H}}^n}\backslash {B_r}} {\int_{{B_r}} {\frac{{v{{\left( \xi  \right)}^{p - 1}}v{{\left( \xi  \right)}^{1 - q}}\psi {{\left( \xi  \right)}^p}}}{{\left\| {{\eta ^{ - 1}} \circ \xi } \right\|_{{{\rm{H}}^n}}^{Q + sp}}}d\xi d\eta } }  + \int_{{{\rm{H}}^n}\backslash {B_r}} {\int_{{B_r}} {\frac{{\left( {u\left( \eta  \right)} \right)_ - ^{p - 1}v{{\left( \xi  \right)}^{1 - q}}\psi {{\left( \xi  \right)}^p}}}{{\left\| {{\eta ^{ - 1}} \circ \xi } \right\|_{{{\rm{H}}^n}}^{Q + sp}}}d\xi d\eta } } } \right)\\
    &= c\left( p \right)\left( {\int_{{{\rm{H}}^n}\backslash {B_r}} {\int_{{B_r}} {\frac{{v{{\left( \xi  \right)}^{p - q}}\psi {{\left( \xi  \right)}^p}}}{{\left\| {{\eta ^{ - 1}} \circ \xi } \right\|_{{{\rm{H}}^n}}^{Q + sp}}}d\xi d\eta } }  + \int_{{{\rm{H}}^n}\backslash {B_R}} {\int_{{B_r}} {\frac{{\left( {u\left( \eta  \right)} \right)_ - ^{p - 1}v{{\left( \xi  \right)}^{1 - q}}\psi {{\left( \xi  \right)}^p}}}{{\left\| {{\eta ^{ - 1}} \circ \xi } \right\|_{{{\rm{H}}^n}}^{Q + sp}}}d\xi d\eta } } } \right).
\end{align*}
Note
\[\omega {\left( \xi  \right)^p} = v{\left( \xi  \right)^{p - q}} = v{\left( \xi  \right)^{p - 1}} \cdot v{\left( \xi  \right)^{1 - q}} \ge {d^{p - 1}} \cdot v{\left( \xi  \right)^{1 - q}},\]
thus
\begin{align}\label{eq64}
   {I_3} \le& c\left( p \right)\left( {\int_{{{\rm{H}}^n}\backslash {B_r}} {\int_{{B_r}} {\frac{{\omega {{\left( \xi  \right)}^p}\psi {{\left( \xi  \right)}^p}}}{{\left\| {{\eta ^{ - 1}} \circ \xi } \right\|_{{{\rm{H}}^n}}^{Q + sp}}}d\xi d\eta } }  + {d^{1 - p}}\int_{{{\rm{H}}^n}\backslash {B_R}} {\int_{{B_r}} {\frac{{\left( {u\left( \eta  \right)} \right)_ - ^{p - 1}\omega {{\left( \xi  \right)}^p}\psi {{\left( \xi  \right)}^p}}}{{\left\| {{\eta ^{ - 1}} \circ \xi } \right\|_{{{\rm{H}}^n}}^{Q + sp}}}d\xi d\eta } } } \right) \nonumber\\
   \le& c\left( {n,p,s} \right)\left( {\mathop {{\rm{ess}}\sup }\limits_{\xi  \in {\rm{supp}}\;\psi } \int_{{{\rm{H}}^n}\backslash {B_r}} {\frac{1}{{\left\| {{\eta ^{ - 1}} \circ \xi } \right\|_{{{\rm{H}}^n}}^{Q + sp}}}d\eta }  + {d^{1 - p}}\int_{{{\rm{H}}^n}\backslash {B_R}} {\frac{{\left( {u\left( \eta  \right)} \right)_ - ^{p - 1}}}{{\left\| {{\eta ^{ - 1}} \circ {\xi _0}} \right\|_{{{\rm{H}}^n}}^{Q + sp}}}d\eta } } \right)\cdot\nonumber\\
   &\int_{{B_r}} {\omega {{\left( \xi  \right)}^p}\psi {{\left( \xi  \right)}^p}d\xi } ,
\end{align}
where we used in the last step that
\[\frac{{{{\left\| {{\eta ^{ - 1}} \circ {\xi _0}} \right\|}_{{{\rm{H}}^n}}}}}{{{{\left\| {{\eta ^{ - 1}} \circ \xi } \right\|}_{{{\rm{H}}^n}}}}} \le c\frac{{{{\left\| {{\eta ^{ - 1}} \circ \xi } \right\|}_{{{\rm{H}}^n}}} + {{\left\| {{\xi ^{ - 1}} \circ {\xi _0}} \right\|}_{{{\rm{H}}^n}}}}}{{{{\left\| {{\eta ^{ - 1}} \circ \xi } \right\|}_{{{\rm{H}}^n}}}}} \le 2c.\]

Now we estimate $I_2$. Without loss of generality, we assume $v\left( \xi  \right) > v\left( \eta  \right)$. By Lemma \ref{Le61}, for any $\varepsilon  \in \left( {0,1} \right)$, it follows
\[\psi {\left( \xi  \right)^p} \le \left( {1 + {c_p}\varepsilon } \right)\psi {\left( \eta  \right)^p} + \left( {1 + {c_p}\varepsilon } \right){\varepsilon ^{1 - p}}{\left| {\psi \left( \xi  \right) - \psi \left( \eta  \right)} \right|^p},\]
where ${c_p} = \left( {p - 1} \right)\Gamma \left( {\max \{ 1,p - 2\} } \right)$. If we take
\begin{center}
$\varepsilon : = \delta \frac{{v\left( \xi  \right) - v\left( \eta  \right)}}{{v\left( \xi  \right)}} \in \left( {0,1} \right),$ for any $\delta  \in \left( {0,1} \right)$,
\end{center}
then
\begin{align}\label{eq65}
   & \frac{{{{\left| {v\left( \xi  \right) - v\left( \eta  \right)} \right|}^{p - 2}}\left( {v\left( \xi  \right) - v\left( \eta  \right)} \right)\left( {v{{\left( \xi  \right)}^{1 - q}}\psi {{\left( \xi  \right)}^p} - v{{\left( \eta  \right)}^{1 - q}}\psi {{\left( \eta  \right)}^p}} \right)}}{{\left\| {{\eta ^{ - 1}} \circ \xi } \right\|_{{{\rm{H}}^n}}^{Q + sp}}}\nonumber \\
  = &  \frac{{{{\left( {v\left( \xi  \right) - v\left( \eta  \right)} \right)}^{p - 1}}}}{{\left\| {{\eta ^{ - 1}} \circ \xi } \right\|_{{{\rm{H}}^n}}^{Q + sp}}}v{\left( \xi  \right)^{1 - q}}\left( {\psi {{\left( \xi  \right)}^p} - \frac{{v{{\left( \eta  \right)}^{1 - q}}\psi {{\left( \eta  \right)}^p}}}{{v{{\left( \xi  \right)}^{1 - q}}}}} \right)\nonumber \\
   \le& \frac{{{{\left( {v\left( \xi  \right) - v\left( \eta  \right)} \right)}^{p - 1}}}}{{\left\| {{\eta ^{ - 1}} \circ \xi } \right\|_{{{\rm{H}}^n}}^{Q + sp}}}v{\left( \xi  \right)^{1 - q}}\psi {\left( \eta  \right)^p}\left( {\left( {1 + {c_p}\varepsilon } \right) - \frac{{v{{\left( \eta  \right)}^{1 - q}}}}{{v{{\left( \xi  \right)}^{1 - q}}}}} \right)\nonumber \\
  & + \left( {1 + {c_p}\varepsilon } \right){\varepsilon ^{1 - p}}\frac{{{{\left( {v\left( \xi  \right) - v\left( \eta  \right)} \right)}^{p - 1}}}}{{\left\| {{\eta ^{ - 1}} \circ \xi } \right\|_{{{\rm{H}}^n}}^{Q + sp}}}v{\left( \xi  \right)^{1 - q}}{\left| {\psi \left( \xi  \right) - \psi \left( \eta  \right)} \right|^p}\nonumber \\
   \le& \frac{{{{\left( {v\left( \xi  \right) - v\left( \eta  \right)} \right)}^{p - 1}}}}{{\left\| {{\eta ^{ - 1}} \circ \xi } \right\|_{{{\rm{H}}^n}}^{Q + sp}}}v{\left( \xi  \right)^{1 - q}}\psi {\left( \eta  \right)^p}\left( {1 + {c_p}\delta \frac{{v\left( \xi  \right) - v\left( \eta  \right)}}{{v\left( \xi  \right)}} - \frac{{v{{\left( \eta  \right)}^{1 - q}}}}{{v{{\left( \xi  \right)}^{1 - q}}}}} \right)\nonumber \\
  & + c\left( p \right){\delta ^{1 - p}}\frac{{v{{\left( \xi  \right)}^{p - q}}{{\left| {\psi \left( \xi  \right) - \psi \left( \eta  \right)} \right|}^p}}}{{\left\| {{\eta ^{ - 1}} \circ \xi } \right\|_{{{\rm{H}}^n}}^{Q + sp}}}.
\end{align}
Noting
\begin{align*}
   1 + {c_p}\delta \frac{{v\left( \xi  \right) - v\left( \eta  \right)}}{{v\left( \xi  \right)}} - \frac{{v{{\left( \eta  \right)}^{1 - q}}}}{{v{{\left( \xi  \right)}^{1 - q}}}} & = {c_p}\delta \frac{{v\left( \xi  \right) - v\left( \eta  \right)}}{{v\left( \xi  \right)}} + \frac{{v{{\left( \xi  \right)}^{1 - q}} - v{{\left( \eta  \right)}^{1 - q}}}}{{v{{\left( \xi  \right)}^{1 - q}}}} \\
   &  = \frac{{v\left( \xi  \right) - v\left( \eta  \right)}}{{v\left( \xi  \right)}}\left( {{c_p}\delta  + \frac{{v\left( \xi  \right)\left( {v{{\left( \xi  \right)}^{1 - q}} - v{{\left( \eta  \right)}^{1 - q}}} \right)}}{{v{{\left( \xi  \right)}^{1 - q}}\left( {v\left( \xi  \right) - v\left( \eta  \right)} \right)}}} \right)\\
   & = \frac{{v\left( \xi  \right) - v\left( \eta  \right)}}{{v\left( \xi  \right)}}\left( {{c_p}\delta  + \frac{{1 - {{\left( {\frac{{v\left( \eta  \right)}}{{v\left( \xi  \right)}}} \right)}^{1 - q}}}}{{\left( {1 - \frac{{v\left( \eta  \right)}}{{v\left( \xi  \right)}}} \right)}}} \right).
\end{align*}
Thus the first term on the right-hand side of \eqref{eq65} satisfies
\begin{align}\label{eq66}
   & \frac{{{{\left( {v\left( \xi  \right) - v\left( \eta  \right)} \right)}^{p - 1}}}}{{\left\| {{\eta ^{ - 1}} \circ \xi } \right\|_{{{\rm{H}}^n}}^{Q + sp}}}v{\left( \xi  \right)^{1 - q}}\psi {\left( \eta  \right)^p}\left( {1 + {c_p}\delta \frac{{v\left( \xi  \right) - v\left( \eta  \right)}}{{v\left( \xi  \right)}} - \frac{{v{{\left( \eta  \right)}^{1 - q}}}}{{v{{\left( \xi  \right)}^{1 - q}}}}} \right) \nonumber\\
   =& \frac{{{{\left( {v\left( \xi  \right) - v\left( \eta  \right)} \right)}^p}}}{{\left\| {{\eta ^{ - 1}} \circ \xi } \right\|_{{{\rm{H}}^n}}^{Q + sp}}}v{\left( \xi  \right)^{ - q}}\psi {\left( \eta  \right)^p}\left( {{c_p}\delta  + \frac{{1 - {{\left( {\frac{{v\left( \eta  \right)}}{{v\left( \xi  \right)}}} \right)}^{1 - q}}}}{{\left( {1 - \frac{{v\left( \eta  \right)}}{{v\left( \xi  \right)}}} \right)}}} \right) = :{J_1}.
\end{align}
For any $t \in \left( {0,1} \right)$, define the real-valued function
\[g\left( t \right): = \frac{{1 - {t^{1 - q}}}}{{1 - t}} =  - \frac{{q - 1}}{{1 - t}}\int\limits_t^1 {{\tau ^{ - q}}d\tau } .\]
Since $q > 1$, we have $g\left( t \right) \le  - \left( {q - 1} \right)$ for all $t \in \left( {0,1} \right)$. If $t \in \left( {0,\frac{1}{2}} \right]$, then
\[g\left( t \right) \le  - \frac{{q - 1}}{{{2^q}}}\frac{{{t^{1 - q}}}}{{1 - t}}.\]

Let us treat two cases: (i) $v\left( \xi  \right) \ge 2v\left( \eta  \right)$ and (ii) $v\left( \eta  \right) < v\left( \xi  \right) < 2v\left( \eta  \right)$.

(i) Let $v\left( \xi  \right) \ge 2v\left( \eta  \right)$, we take $t = \frac{{v\left( \eta  \right)}}{{v\left( \xi  \right)}}$. Then $t \in \left( {0,\frac{1}{2}} \right]$ and thus
\begin{align}\label{eq67}
   {J_1}& \le \frac{{{{\left( {v\left( \xi  \right) - v\left( \eta  \right)} \right)}^p}\psi {{\left( \eta  \right)}^p}}}{{\left\| {{\eta ^{ - 1}} \circ \xi } \right\|_{{{\rm{H}}^n}}^{Q + sp}v{{\left( \xi  \right)}^q}}}\left( {{c_p}\delta  - \frac{{q - 1}}{{{2^q}}}\frac{{v{{\left( \xi  \right)}^{q - 1}}}}{{v{{\left( \eta  \right)}^{q - 1}}}}\frac{{v\left( \xi  \right)}}{{v\left( \xi  \right) - v\left( \eta  \right)}}} \right)\nonumber \\
   &  \le \frac{{{{\left( {v\left( \xi  \right) - v\left( \eta  \right)} \right)}^{p - 1}}\psi {{\left( \eta  \right)}^p}}}{{\left\| {{\eta ^{ - 1}} \circ \xi } \right\|_{{{\rm{H}}^n}}^{Q + sp}v{{\left( \eta  \right)}^{q - 1}}}}\left( {{c_p}\delta  - \frac{{q - 1}}{{{2^q}}}} \right),
\end{align}
where we used $\frac{{v{{\left( \eta  \right)}^{q - 1}}\left( {v\left( \xi  \right) - v\left( \eta  \right)} \right)}}{{v{{\left( \xi  \right)}^q}}} \le 1$ in the last step. Taking $\delta  = \frac{{q - 1}}{{{c_p}{2^{q + 1}}}}$, we have from \eqref{eq67} that
\begin{equation}\label{eq68}
{J_1} \le  - \frac{{q - 1}}{{{2^{q + 1}}}}\frac{{{{\left( {v\left( \xi  \right) - v\left( \eta  \right)} \right)}^{p - 1}}\psi {{\left( \eta  \right)}^p}}}{{\left\| {{\eta ^{ - 1}} \circ \xi } \right\|_{{{\rm{H}}^n}}^{Q + sp}v{{\left( \eta  \right)}^{q - 1}}}}.
\end{equation}
Moreover, since $v\left( \xi  \right) \ge 2v\left( \eta  \right)$, we see
\begin{align*}
   \frac{{{{\left( {v\left( \xi  \right) - v\left( \eta  \right)} \right)}^{p - 1}}}}{{v{{\left( \eta  \right)}^{q - 1}}}}& \ge \frac{{{2^{q - 1}}{{\left( {v\left( \xi  \right) - v\left( \eta  \right)} \right)}^{p - 1}}}}{{v{{\left( \xi  \right)}^{q - 1}}}} \ge {2^{q - p}}v{\left( \xi  \right)^{p - q}} \\
   &  \ge {2^{q - p}}{\left( {v{{\left( \xi  \right)}^{\frac{{p - q}}{p}}} - v{{\left( \eta  \right)}^{\frac{{p - q}}{p}}}} \right)^p},
\end{align*}
and obtain by combining with \eqref{eq68} that
\begin{equation}\label{eq69}
{J_1} \le  - \frac{{q - 1}}{{{2^{p + 1}}}}\frac{{{{\left( {\omega \left( \xi  \right) - \omega \left( \eta  \right)} \right)}^p}\psi {{\left( \eta  \right)}^p}}}{{\left\| {{\eta ^{ - 1}} \circ \xi } \right\|_{{{\rm{H}}^n}}^{Q + sp}}}.
\end{equation}

(ii) Let $v\left( \eta  \right) < v\left( \xi  \right) < 2v\left( \eta  \right)$, and take $t = \frac{{v\left( \eta  \right)}}{{v\left( \xi  \right)}}$, $\delta  = \frac{{q - 1}}{{{c_p}{2^{q + 1}}}}$, then $t \in \left( {\frac{1}{2},1} \right)$ and
\begin{align}\label{eq610}
   {J_1}& \le \frac{{{{\left( {v\left( \xi  \right) - v\left( \eta  \right)} \right)}^p}\psi {{\left( \eta  \right)}^p}}}{{\left\| {{\eta ^{ - 1}} \circ \xi } \right\|_{{{\rm{H}}^n}}^{Q + sp}v{{\left( \xi  \right)}^q}}}\left( {{c_p}\delta  - \left( {q - 1} \right)} \right)\nonumber \\
   &  =  - \left( {{2^{q + 1}} - 1} \right)\frac{{q - 1}}{{{2^{q + 1}}}}\frac{{{{\left( {v\left( \xi  \right) - v\left( \eta  \right)} \right)}^p}\psi {{\left( \eta  \right)}^p}}}{{\left\| {{\eta ^{ - 1}} \circ \xi } \right\|_{{{\rm{H}}^n}}^{Q + sp}v{{\left( \xi  \right)}^q}}}.
\end{align}
Noting
\begin{align*}
   {\left( {\omega \left( \xi  \right) - \omega \left( \eta  \right)} \right)^p}& = {\left( {\frac{{p - q}}{p}} \right)^p}{\left( {\int\limits_{v\left( \eta  \right)}^{v\left( \xi  \right)} {{t^{ - \frac{q}{p}}}dt} } \right)^p} \le {\left( {\frac{{p - q}}{p}} \right)^p}\frac{1}{{v{{\left( \eta  \right)}^q}}}{\left( {\int\limits_{v\left( \eta  \right)}^{v\left( \xi  \right)} {dt} } \right)^p} \\
   &  = {\left( {\frac{{p - q}}{p}} \right)^p}\frac{{{{\left( {v\left( \xi  \right) - v\left( \eta  \right)} \right)}^p}}}{{v{{\left( \eta  \right)}^q}}} \le {2^q}{\left( {\frac{{p - q}}{p}} \right)^p}\frac{{{{\left( {v\left( \xi  \right) - v\left( \eta  \right)} \right)}^p}}}{{v{{\left( \xi  \right)}^q}}},
\end{align*}
we deduce by combining \eqref{eq610} that
\begin{equation}\label{eq611}
{J_1} \le  - \left( {{2^{q + 1}} - 1} \right)\frac{{q - 1}}{{{2^{2q + 1}}}}{\left( {\frac{p}{{p - q}}} \right)^p}\frac{{{{\left( {\omega \left( \xi  \right) - \omega \left( \eta  \right)} \right)}^p}\psi {{\left( \eta  \right)}^p}}}{{\left\| {{\eta ^{ - 1}} \circ \xi } \right\|_{{{\rm{H}}^n}}^{Q + sp}}}.
\end{equation}

In conclusion, when $v\left( \xi  \right) > v\left( \eta  \right)$, it yields by comparing \eqref{eq69} and \eqref{eq611} that
\begin{equation}\label{eq612}
{J_1} \le  - c\frac{{{{\left( {\omega \left( \xi  \right) - \omega \left( \eta  \right)} \right)}^p}\psi {{\left( \eta  \right)}^p}}}{{\left\| {{\eta ^{ - 1}} \circ \xi } \right\|_{{{\rm{H}}^n}}^{Q + sp}}},
\end{equation}
where
\[c = \min \left\{ {\frac{{q - 1}}{{{2^{p + 1}}}},\left( {{2^{q + 1}} - 1} \right)\frac{{q - 1}}{{{2^{2q + 1}}}}{{\left( {\frac{p}{{p - q}}} \right)}^p}} \right\}.\]

When $v\left( \xi  \right) = v\left( \eta  \right)$, \eqref{eq612} obviously holds.

When $v\left( \xi  \right) < v\left( \eta  \right)$, the discussion is similar; it suffices to swap the positions of $\xi$ and $\eta$.

Thus, for all cases, we have
\begin{align}\label{eq613}
  {I_2} \le& - c\left( {p,q} \right)\int_{{B_r}} {\int_{{B_r}} {\frac{{{{\left( {\omega \left( \xi  \right) - \omega \left( \eta  \right)} \right)}^p}\psi {{\left( \eta  \right)}^p}}}{{\left\| {{\eta ^{ - 1}} \circ \xi } \right\|_{{{\rm{H}}^n}}^{Q + sp}}}d\xi d\eta } }\nonumber  \\
   &  + \frac{{c\left( p \right)}}{{{{\left( {q - 1} \right)}^{p - 1}}}}\int_{{B_r}} {\int_{{B_r}} {\frac{{{{\left( {\max \left\{ {\omega \left( \xi  \right),\omega \left( \eta  \right)} \right\}} \right)}^p}{{\left| {\psi \left( \xi  \right) - \psi \left( \eta  \right)} \right|}^p}}}{{\left\| {{\eta ^{ - 1}} \circ \xi } \right\|_{{{\rm{H}}^n}}^{Q + sp}}}d\xi d\eta } } .
\end{align}
Combining \eqref{eq62}-\eqref{eq64} and \eqref{eq613}, we conclude \eqref{eq61}.
\end{proof}

\textbf{Proof of Theorem \ref{Th14}.} Let $r \in \left( {0,1} \right],\;\frac{1}{2} < \tau ' < \tau  \le \frac{3}{4}$, and choose a cut-off function $\psi  \in C_0^\infty \left( {{B_{\tau r}}\left( {{\xi _0}} \right)} \right)$ satisfying $0 \le \psi  \le 1,\;\left| {{\nabla _H}\psi } \right| \le \frac{4}{{\left( {\tau  - \tau '} \right)r}}$ in ${B_{\tau r}}\left( {{\xi _0}} \right)$ and $\psi  = 1$ in ${B_{\tau 'r}}\left( {{\xi _0}} \right)$. For $1 < q < p,\;d > 0$, denote
\[v = u + d,\;\omega  = {\left( {u + d} \right)^{\frac{{p - q}}{p}}}.\]
Using the properties of $\psi $, we obtain
\begin{equation}\label{eq614}
{I_1}: = \int_{{B_r}} {{\omega ^p}{{\left| {{\nabla _H}\psi } \right|}^p}d\xi }  \le \frac{{c{r^{ - p}}}}{{{{\left( {\tau  - \tau '} \right)}^p}}}\int_{{B_{\tau r}}} {{\omega ^p}d\xi } .
\end{equation}
By the estimate in (4.28) of \cite{PP22}, we have
\begin{align}\label{eq615}
   {I_2}&: = \int_{{B_r}} {\int_{{B_r}} {\frac{{{{\left( {\max \left\{ {\omega \left( \xi  \right),\omega \left( \eta  \right)} \right\}} \right)}^p}{{\left| {\psi \left( \xi  \right) - \psi \left( \eta  \right)} \right|}^p}}}{{\left\| {{\eta ^{ - 1}} \circ \xi } \right\|_{{{\rm{H}}^n}}^{Q + sp}}}d\xi d\eta } } \nonumber \\
   &  \le \frac{{c{r^{ - sp}}}}{{{{\left( {\tau  - \tau '} \right)}^p}}}\int_{{B_{\tau r}}} {{\omega ^p}d\xi }  \le \frac{{c{r^{ - p}}}}{{{{\left( {\tau  - \tau '} \right)}^p}}}\int_{{B_{\tau r}}} {{\omega ^p}d\xi } .
\end{align}

Let ${\rm{Tail}}({u_ - };{\xi _0},R) > 0$, then for any $\varepsilon  > 0$ and $r \in \left( {0,1} \right]$, take
\[d = \frac{1}{2}{\left( {\frac{r}{R}} \right)^{\frac{p}{{p - 1}}}}{\rm{Tail}}({u_ - };{\xi _0},R) + \varepsilon  > 0.\]
Note
\begin{equation}\label{eq616}
\mathop {{\rm{ess}}\sup }\limits_{\xi  \in {\rm{supp}}\;\psi } \int_{{{\rm{H}}^n}\backslash {B_r}} {\frac{1}{{\left\| {{\eta ^{ - 1}} \circ \xi } \right\|_{{{\rm{H}}^n}}^{Q + sp}}}d\eta }  \le c\left( {n,p,s} \right){r^{ - sp}} \le c\left( {n,p,s} \right){r^{ - p}},
\end{equation}
so
\begin{align}\label{eq617}
   {I_3}&: = \left( {\mathop {{\rm{ess}}\sup }\limits_{\xi  \in {\rm{supp}}\;\psi } \int_{{{\rm{H}}^n}\backslash {B_r}} {\frac{1}{{\left\| {{\eta ^{ - 1}} \circ \xi } \right\|_{{{\rm{H}}^n}}^{Q + sp}}}d\eta }  + {d^{1 - p}}{R^{ - p}}{\rm{Tail}}{{({u_ - };{\xi _0},R)}^{p - 1}}} \right)\int_{{B_r}} {\omega {{\left( \xi  \right)}^p}\psi {{\left( \xi  \right)}^p}d\xi } \nonumber \\
   &  \le \frac{{c\left( {n,p,s} \right){r^{ - p}}}}{{{{\left( {\tau  - \tau '} \right)}^p}}}\int_{{B_{\tau r}}} {{\omega ^p}d\xi } .
\end{align}
If ${\rm{Tail}}({u_ - };{\xi _0},R) = 0$, then we can take $d = \varepsilon  > 0$ and use \eqref{eq616} to obtain \eqref{eq617}.

Using $\psi  = 1$ in ${B_{\tau 'r}}\left( {{\xi _0}} \right)$, Lemma \ref{Le22}, Lemma \ref{Le62}, \eqref{eq614}, \eqref{eq615} and \eqref{eq617}, we get
\begin{align}\label{eq618}
   & {\left( {\fint_{{B_{\tau 'r}}} {{\nu ^{\frac{{Q\left( {p - q} \right)}}{{Q - p}}}}d\xi } } \right)^{\frac{{Q - p}}{Q}}} = {\left( {\fint_{{B_{\tau 'r}}} {{\omega ^{\frac{{Qp}}{{Q - p}}}}d\xi } } \right)^{\frac{{Q - p}}{Q}}} \nonumber\\
\le   &  {\left( {\fint_{{B_{\tau r}}} {{{\left| {\omega \psi } \right|}^{\frac{{Qp}}{{Q - p}}}}d\xi } } \right)^{\frac{{Q - p}}{Q}}} \le c{\left( {\tau r} \right)^p}\fint_{{B_{\tau r}}} {{{\left| {{\nabla _H}\left( {\omega \psi } \right)} \right|}^p}d\xi }\nonumber\\
 \le& c{\left( {\tau r} \right)^p}\left( {\fint_{{B_{\tau r}}} {{{\left| {{\nabla _H}\omega } \right|}^p}{\psi ^p}d\xi }  + \fint_{{B_{\tau r}}} {{{\left| {{\nabla _H}\psi } \right|}^p}{\omega ^p}d\xi } } \right)\nonumber\\
 \le& \frac{c}{{{{\left( {\tau  - \tau '} \right)}^p}}}\fint_{{B_{\tau r}}} {{\omega ^p}d\xi } ,
\end{align}
where $c = c\left( {n,p,s} \right)$. Using \eqref{eq618}, $1 < q < p$, and the standard Moser iteration technique, it yields
\begin{equation}\label{eq619}
{\left( {\fint_{{B_{\frac{r}{2}}}} {{\nu ^l}d\xi } } \right)^{\frac{1}{l}}} \le c{\left( {\fint_{\frac{{3r}}{4}} {{\nu ^{l'}}d\xi } } \right)^{\frac{1}{{l'}}}},\;0 < l' < l < \frac{{Q\left( {p - 1} \right)}}{{Q - p}}.
\end{equation}
Let $\lambda  \in \left( {0,1} \right)$ be as in Lemma \ref{Le53}, take $l' = \lambda  \in \left( {0,1} \right)$ and use
\[{\left( {\fint_{{B_{\frac{r}{2}}}} {{u^l}d\xi } } \right)^{\frac{1}{l}}} \le {\left( {\fint_{{B_{\frac{r}{2}}}} {{\nu ^l}d\xi } } \right)^{\frac{1}{l}}},\]
\eqref{eq619} and Lemma \ref{Le53}, we have that for all $0 < l < \frac{{Q\left( {p - 1} \right)}}{{Q - p}}$,
\begin{equation}\label{eq620}
 {\left( {\fint_{{B_{\frac{r}{2}}}} {{u^l}d\xi } } \right)^{\frac{1}{l}}} \le c\mathop {{\rm{ess inf}}}\limits_{{B_r}\left( {{\xi _0}} \right)} \nu  + c{\left( {\frac{r}{R}} \right)^{\frac{p}{{p - 1}}}}{\rm{Tail}}({u_ - };{\xi _0},R),
\end{equation}
where $v = u + d$. For any $\varepsilon  > 0$ and $r \in \left( {0,1} \right]$, take
\[d = \frac{1}{2}{\left( {\frac{r}{R}} \right)^{\frac{p}{{p - 1}}}}{\rm{Tail}}({u_ - };{\xi _0},R) + \varepsilon \]
in \eqref{eq620}. Then letting $\varepsilon  \to 0$ yields \eqref{eq19}.

\section*{Acknowledgements}
This work was supported by the National Natural Science Foundation of China (No. 12501269) and the Natural Science Basic Research Program of Shaanxi (No. 2024JC-YBQN-0054).

\section*{Declarations}
\subsection*{Conflict of interest} The authors declare that there is no conflict of interest. We also declare that this
manuscript has no associated data.

\subsection*{Data Availability} Data sharing is not applicable to this article as no datasets were generated or analysed during the current study.


\end{document}